\documentclass[review]{elsarticle}

\usepackage{graphicx}%
\usepackage{multirow}%
\usepackage{amsmath,amssymb,amsfonts}%
\usepackage{amsthm}%
\usepackage{mathrsfs}%
\usepackage[title]{appendix}%
\usepackage{xcolor}%
\usepackage{textcomp}%
\usepackage{manyfoot}%
\usepackage{booktabs}%
\usepackage{algorithm}%
\usepackage{algorithmicx}%
\usepackage{algpseudocode}%
\usepackage{listings}%

\usepackage{hyperref}
\usepackage{amsmath}
\usepackage[numbers]{natbib} 
\usepackage[latin1]{inputenc}
\usepackage[english]{babel}
\usepackage{epsfig}
\usepackage{color}
\usepackage{amssymb}
\usepackage{float}
\usepackage{float}
\usepackage{subfigure}
\usepackage{graphics}
\usepackage{graphicx}
\usepackage{indentfirst}
\usepackage{amsthm}
\usepackage{mathabx}
\usepackage[percent]{overpic}
\usepackage{tikz}

\newtheorem{theorem}{Theorem}[section]
\newtheorem{definition}{Definition}[section]
\newtheorem{example}{Example}
\newtheorem{proposition}{Proposition}[section]

\newtheorem{remark}{Remark}

\newcommand{\R}{\mathbb{R}}
\newcommand{\C}{\mathbb{C}}
\newcommand{\RF}{\mathbb{R}_\mathcal{F}}

\newcommand{\RFA}{\mathbb{R}_{\mathcal{F}(A)}}

\journal{Journal of \LaTeX\ Templates}

\begin{document}

\begin{frontmatter}

\title{Calculus for Functions with Fuzzy Inputs and Outputs: Applications to Fuzzy Differential Equations}
%\tnotetext[mytitlenote]{Fully documented templates are available in the elsarticle package on \href{http://www.ctan.org/tex-archive/macros/latex/contrib/elsarticle}{CTAN}.}

%% Group authors per affiliation:
\author{La\'ecio Carvalho de Barros\fnref{myfootnote2}}
\address{University of Campinas, Department of Applied Mathematics, 13081-970, Campinas, SP, Brazil}
\fntext[myfootnote2]{CNPq grant no. 314885/2021-8}

\author{Estev\~ao Esmi \fnref{myfootnote}}
\address{University of Campinas, Department of Applied Mathematics, 13081-970, Campinas, SP, Brazil}
\fntext[myfootnote]{CNPq grant no. 311976/2023-9 and Fapesp grant no. 2020/09838-0}

\author{Francielle Santo Pedro Sim\~oes}%\fnref{myfootnote1}}
\address{Federal University of S\~ao Paulo, Multidisciplinary Department, 06110-295, Osasco, SP, Brazil}
%\fntext[myfootnote1]{}

\author{Mina Shahidi \fnref{myfootnote3}}
\address{University of Campinas, Department of Applied Mathematics, 13081-970, Campinas, SP, Brazil}
\fntext[myfootnote3]{Fapesp grant no. 2022/00196-1}

\begin{abstract}
This article presents a theory of differential and integral calculus for mapping between Banach spaces formed by subsets of fuzzy numbers called $A$-linearly correlated fuzzy numbers ($\RFA$), where both the domain and codomain are spaces composed of fuzzy numbers. This is one of the main contributions of this study from a theoretical point of view, as well-known approaches to fuzzy calculus in the literature typically deal with fuzzy number-valued functions defined on intervals of real numbers. Notions of differentiability and integrability based on complex functions are proposed. Moreover, we introduce the study of ordinary differential equations for which the solutions are functions $\RFA \to \RFA$ or $D \subseteq \R \to \RFA$. For the latter case, we present an initial study of the solution and its phase portrait for two-dimensional differential equation systems. In particular, for the former case, we examine the Lotka-Volterra model and analyze its phase portrait.
\end{abstract}

\begin{keyword}
Calculus for fuzzy mapping, Cauchy-Riemann equations, Fuzzy chain rule, Fuzzy analytic functions, %Fuzzy curve, 
Fuzzy differential equation, Fuzzy Lotka-Volterra
\end{keyword}

\end{frontmatter}

%\linenumbers

\section{Introduction}\label{sec1}

Most approaches to fuzzy calculus presented in the literature are developed for functions whose domain is an interval of $\R$ and the codomain is the set of fuzzy numbers, denoted here by $\RF$, 
(e.g. see \cite{de2017fuzzy, bede13, esmi2023some, gomes2015fuzzy,Kaleva87, Kaleva06, lakshmikantham2004theory, mazandarani2017granular, Puri, Seikkala87, Seikkala02} and the references therein). In contrast, this work deals with differential and integral calculus theory for mappings between two certain subclasses of fuzzy numbers. Hence, to the best of our knowledge, this is the first time that a calculus theory is proposed for %fuzzy number-valued functions defined in a subset of fuzzy numbers, that is, the case 
functions where both the input and output %of the function 
are fuzzy numbers. 
To this end, we shall consider the space $\RFA$ introduced by Esmi {\it et al.} \cite{de2021differential, esmi2018frechet, santo2020calculus} which is a Banach space isomorphic to $\R^2$ such that the vectors are elements of $\RF$ that are linearly interactive with each other.
%with a given fuzzy number $A$. 
Using the $\RFA$ space, many works were developed \cite{pedro2023fuzzy, shen2022calculus, shen2022first, simoes2023interactive, wasques2020interactive, wasques2018higher} including a generalization of this space to higher orders \cite{esmi2021banach,  esmi2022calculus, laiate2021bidimensional, laiate2021cross}. 
However, all these articles only deal with fuzzy number-valued functions defined in an interval of real numbers. 
In this approach, the term fuzzy curve refers to fuzzy number-valued functions defined on an interval of $\R$, and the term fuzzy mapping refers to functions that map fuzzy numbers to fuzzy numbers.

Problems where the input may be fuzzy for a variety of reasons, such as intrinsic uncertainties associated with assessing the value of the independent variable, among others, arise frequently in mathematical modeling.

For example, the time variable in the study of the spread of diseases or wildfires can also be considered fuzzy, as the moment of infection or the initiation of a fire is uncertain. Another example is the study of the concentration of substances at the border (or at the bottom) of a river, which also may be modeled as a phenomenon with fuzzy inputs since the border (or the bottom) of a river is uncertain and can be described by fuzzy sets. 

From a mathematical point of view, a motivation for studying functions that map fuzzy numbers to fuzzy numbers is that it is possible to confront compositions $f \circ g$, where $g: \mathbb{R} \to \mathbb{\RF}$ and $f: \mathbb{\RF} \to \mathbb{\RF}$. 
Note that $f \circ g$ is a function from $\mathbb{R}$ to $\mathbb{\RF}$, and a natural mathematical question that arises is: Is it possible to establish, in some sense, a type of chain rule for the derivative of $f \circ g$?
It is worth noting that compositions of functions naturally appear in the study of differential equations, as the $n$th derivative of a fuzzy number-valued function $x$ is described by the composition of $x$ and its derivatives of lower order with a given function $f$. A common method to study the stability of equilibrium points  (i.e, the points $\bar{x}$ such that $f(\bar{x}) = 0$) in a first-order system of classical (crisp) differential equations involves assessing the derivative of $f$ with respect to $x$.

Beyond the fact that the formulation of many problems may involve products of fuzzy numbers, the notions of multiplication operation and reciprocal are also important if one intends to define derivative as a limit of quotients. 
In this context, the usual multiplication operator on the class of fuzzy numbers, which coincides with the Zadeh extension of the real number multiplication operator, can not be used here because $\RFA$ is not closed with respect to it. 
In other words, for $B,C \in \RFA$, it may happen that the usual product of $B$ and $C$ does not belong to $\RFA$. 
Recently, this problem has been overpassed by using the $\Psi$-cross product operator which is defined on $\RFA$ whenever $1$-level of  $A$, denoted by $[A]_1$, is a unit set, say $[A]_1 = \{a\}$ \cite{laiate2021cross,longo2022cross}.
In addition, the $\Psi$-cross product can be obtained in terms of the extension principle of the linearization of the multiplication of real numbers around the pair of points for which the operands reach the maximum membership \cite{laiate2021cross}. However, the notion of a reciprocal or inverse element, that is, for $B \in \RFA$ there exists $B^{-1} \in \RFA$ such that the $\Psi$-cross product of $B$ and $B^{-1}$ is equal to the real number $1$, is not defined for all $B \in \RFA $ such that $B$ does not represent the number zero ($B\neq 0$).
For example, if $a = 0$, which is a case of great interest from a modeling point of view because $A$ can be interpreted as fuzzy noise around zero, then for every $B = qA$, $q \in \R$ does not have a reciprocal \cite{laiate2021cross}.
In contrast, this paper introduces the multiplication $\odot$, which is induced by the multiplication of complex numbers (see \eqref{multiplicacao}) inheriting all its properties such as every element of $\RFA \setminus \{0\}$ has a reciprocal, it extends the multiplication on $\R$, etc. Therefore, the proposed multiplication can be used to define derivative in terms of the limit of quotients as aforementioned.

The article is organized as follows. Section \ref{Premilinary} presents some basic concepts and notations of fuzzy set theory. Section \ref{Thespace} presents the space of A-Linearly Correlated Fuzzy Numbers and proves that it is isomorphic to the field $\C$.  Section \ref{Analyticalfuzzyfunctions}  details a theory of differential and integral calculus for analytical functions in $\RFA$. 
Section \ref{Fuzzyfunctions} introduces the study of fuzzy curves from $\R$ to $\RFA$ and of ordinary differential equations.  In Section 6, we show some applications in two-dimensional differential equation systems
and, finally, in Section 7, we include some conclusions.

\section{Mathematical Background}\label{Premilinary}

Let $\RF$ be the set of fuzzy numbers. For the sake of simplicity, we use the same symbol to denote a fuzzy number and its membership function, that is, for $A\in \RF$ and for all $x \in \R$, the membership degree of $x$ in $A$ is given by $A(x)$. Recall that the $\alpha$-levels of $A$ are non-empty close intervals of $\R$ \cite{de2017first, wasques2018higher}, that we denote here by $[A]_{\alpha} := [\underline{a}_\alpha, \overline{a}_\alpha]$ for all $\alpha \in [0,1]$. Moreover, the space $\RF$ forms a complete metric space with the metric $d_\infty$ defined for every $A,B \in \RF$ as follows \cite{diamond2000metric}:
$$d_{\infty }(A,B)=\sup_{0\leq \alpha \leq 1}\max \{|a^{-}_{\alpha}-b^{-}_{\alpha}|,|a^{+}_{\alpha}-b^{+}_{\alpha}|\}.$$

Trapezoidal and triangular fuzzy numbers are well-known and used examples of fuzzy numbers that can be defined levelwise as follows. For $a,b,c,d \in \mathbb{R}$ such that $a\leq b\leq c\leq d$, we can defined a trapezoidal fuzzy number $A$, which is also denoted by the symbol $(a;b;c;d)$, such that $[A]_\alpha =[a+\alpha(b-a), d-\alpha(d-c)]$ for all $\alpha \in [0,1]$. 
If $b=c$, we speak of a triangular fuzzy number and denote by the symbol $(a;b;d)$ instead.

The usual addition and scalar product can be defined level-wise as follows:
\begin{equation}\label{somausual}
[C + D]_\alpha = [\underline{c}_\alpha + \underline{d}_\alpha, \overline{c}_\alpha + \overline{d}_\alpha]    
\end{equation}
and 
\begin{equation}\label{prodescalarusual}
[\lambda C]_\alpha = \left\lbrace \begin{array}{cl}
 \left[\lambda \underline{c}_\alpha, \lambda\overline{c}_\alpha\right]    &, \mbox{ if }  \lambda \geq 0, \\
 \left[\lambda \overline{c}_\alpha, \lambda\underline{c}_\alpha\right]    &, \mbox{ if }  \lambda < 0. 
\end{array}\right.     
\end{equation}
for all $B,C \in \RF$ and $\lambda \in \R$.

We say that a fuzzy number $A$ is asymmetric (or nonsymmetric) if for all $x \in \R$ there exists $y \in \R$ such that $A(x-y) \neq A(x+y)$. 
An alternative definition of this concept of asymmetry is provided by Lemma 3.5 in \cite{esmi2018frechet} which states that $A$ is asymmetric if $\underline{a}_\alpha +  \overline{a}_\alpha$ is not constant with respect to $\alpha$. 
An immediate consequence of this lemma is that the set of asymmetric fuzzy numbers is open. Since $A$ is asymmetric, there exist $\alpha \neq \beta$ 
such that $\epsilon:= |(\underline{a}_\alpha + \overline{a}_\alpha) - (\underline{a}_\beta +  \overline{a}_\beta)| > 0$. For all $B \in \RF$ with $d_\infty(B,A) < 0.2\epsilon$ we obtain $|(\underline{b}_\alpha + \overline{b}_\alpha) - (\underline{b}_\beta +  \overline{b}_\beta)| \geq 0.2\epsilon  > 0$. Therefore, all elements in the open ball of radius $0.2\epsilon$ centered at $A$ are also asymmetric.   
One can also prove that the set of asymmetric fuzzy numbers is dense in $\RF$ w.r.t. $d_\infty$ \cite{esmi2021banach}. Thus, from a practical point of view, this set is sufficiently large to approximate any fuzzy number and, furthermore, asymmetry is preserved under small perturbations. 

In \cite{esmi2018frechet}, Esmi et al. prove that the mapping $(r,q)\mapsto r + qA$ is injective if, and only if, $A$ is asymmetric. 
%Let $\RFA$ denote the range of $\Psi$. Note that $\R$ can be embedded in $\RFA$ since $\Psi_A(r,0) = r$ for every real number $r$.   
By composing this mapping and the usual bijection between $\R^2$ and $\C$, we obtain a bijection from $\C$ to its range.  %,that can be used to induce an isometric isomorphism from $\C$ to $\RFA$. 
The following section focuses on this bijection.   

\section{The space of $A$-Linearly Correlated Fuzzy Numbers -$\RFA$}\label{Thespace}
 
 For each $A\in \RF$, we can define the operator $\Phi:\C\to\RFA $ that associates 
each complex number $(r+iq)\in\C$ with the fuzzy number $\Phi_A(r+iq)$ given by 
\begin{equation}\label{op_phi}
\Phi_A(r+iq) = r+qA,  q,r\in\R.
\end{equation}
For notational convenience, we use the symbol 
$\Phi$ instead of $\Phi_A$. 

As in \cite{esmi2018frechet}, it is easy to verify that $ \Phi $ is well defined since $\Psi$ is defined in terms of the usual arithmetic operations. Moreover, one can easily verify that  $[\Phi(r+iq)]_{\alpha} = \{r+qx \in \mathbb{R} \,|\, x \in [A]_{\alpha} \}$ for all $\alpha \in [0,1].$ Let us denote the range of the operator $\Phi$  by the symbol $\mathbb{R}_{\mathcal{F}(A)}$, that is, $\mathbb{R}_{\mathcal{F}(A)} = \{ \Phi(r+iq) \mid (r+iq) \in \mathbb{C} \}$. 

An important fact that is worth emphasizing is that $\mathbb{R}$ can be embedded in $\mathbb{R}_{\mathcal{F}(A)}$ since every real number $ r $ can be identified with the fuzzy number $\Phi(r+i0) \in \mathbb{R}_{\mathcal{F}(A)}$, that is, $\mathbb{R} \subseteq \mathbb{R}_{\mathcal{F}(A)}$. 
Another interesting property is that for every $B \in \mathbb{R}_{\mathcal{F}(A)} \setminus \mathbb{R}$, we have $\mathbb{R}_{\mathcal{F}(B)}= \mathbb{R}_{\mathcal{F}(A)}$ \cite{esmi2018frechet}. In other words, 
the range of $\Psi_B$ does not depend on the choice of $B  \in \mathbb{R}_{\mathcal{F}(A)} \setminus \mathbb{R}$ and is equal to $\RFA$. 

\begin{theorem}\label{teo_injetividade}
Given $A\in \RF$, consider the operator $\Phi:\C \to \RFA$ given by 
\begin{equation}
\Phi(r + iq) =  r + qA    
\end{equation}
for all $z = (r + iq) \in \C$. 
The operator  $\Phi$ is a bijection if, and only if, $A$ is asymmetric.
\end{theorem}
\begin{proof}
It follows analogously from Theorem 3.6 of \cite{esmi2018frechet}.
\end{proof}

Informally and roughly speaking, the space $\RFA$ is obtained from $\C$ simply replacing the imaginary unit $i$ by the asymmetric fuzzy number $A$. However, unlike the complex case, the quantity $r+qA$ can be seen as a sum of fuzzy numbers $r$ and $qA$ that results in a fuzzy number $B=r+qA$. 
Note that, in general, any fuzzy number $B$ can be written as $b + B_0$, where $b$ is a real number and $B_0(x-b) =  B(x)$ for all $x$. For example, the triangular fuzzy number $(0.5;2;5)$ is equal to $2+3(-0.5;0;1)$. However, for a fixed fuzzy number $A$, only elements of $\RFA$ can be decomposed as a sum of a real part $r$ and a fuzzy part $qA$. 

For an asymmetric $A$, if the codomain of the operator $\Phi$ is restricted to its range  ($\mathbb{R}_{\mathcal{F}(A)}$), then $\Phi$ is a bijection from Theorem \ref{teo_injetividade}. 
Using this bijection, one can induce isometric isomorphism from the field $\C$ to $\RFA$. This is stated in the following theorem. 

\begin{theorem}\label{col_ev}
Let $A$ be an asymmetric fuzzy number. The set $\RFA$ is a vector space over the scalar field $\mathbb{R}$ with the usual scalar product and the addition operator given by %respectively by
\begin{eqnarray*}
B \oplus C  =   \Phi(\Phi^{-1}(B) + \Phi^{-1}(C)) =   (r+s) + (q+p)A
\end{eqnarray*}
for all $\eta\in\mathbb{R}$, $B=(r + qA) \in \RFA$ and $C=(s + pA) \in \RFA$.
In addition, $\RFA$ is isometric to $\mathbb{C}$ with the norm 
\begin{equation}
\|B\|_\Phi = |\Phi^{-1}(B)| = \sqrt{r^2+q^2}.%|r| + |p|.      
\end{equation}
\end{theorem}
\begin{proof}
First of all, note that for every $\eta \in \R$ and $B = (r + qA) \in \RFA$, we have that 
\begin{eqnarray*}
 \eta B & = & \eta(r + qA) \\ 
  & = & (\eta r) + (\eta q)A \\
  & = & \Phi((\eta r) + (\eta q)i) \\
  & = & \Phi(\eta( r + iq)) \\
  & = & \Phi(\eta\Phi^{-1}(B)).
\end{eqnarray*}
This implies that the usual scalar product coincides with the scalar product induced by the bijection $\Phi$. 
Therefore, the bijection $\Phi$ becomes an isometric isomorphism from $\C$ to $\RFA$ over the field $\R$ since 
all operations and the norm on $\RFA$ are induced by the bijection $\Phi$.      
\end{proof}

Notice that if $r\in\R \subset \RFA \subset \mathbb{R}_{\mathcal{F}}$, we have that $r\oplus B=r+B$,  %and $r\cdot B=r\cdot B$
 where ``$+$" is the usual addition given by  \eqref{somausual}.
 %and ``$\cdot$" is given by equations \eqref{somausual} and \eqref{prodescalarusual} respectively.
 As usual for every vector space, for every $B,C \in \RFA$, we simply denote the product $(-1)B$ by $-B$ and 
 the difference of $C$ and $B$ is given by
 \begin{equation}\label{eq:subtraction}
 C \ominus B = C \oplus (-B).    
 \end{equation}
 
From now on, unless otherwise stated, we assume that the fuzzy number $A$ that generates the space $\RFA$ is asymmetric and, therefore, $\RFA$ is a vector space with the operations given as in Theorem \ref{col_ev}. 

\begin{remark}\label{rm:division}
Through the isomorphism $\Phi$, we can define the product and division operators between the fuzzy numbers $B=(r + qA) \in \RFA$ and $C=(s + pA) \in \RFA$, as the set of complex numbers $\mathbb{C}$ already possesses these operations \cite{brown2009complex}. More precisely, we have   
\begin{eqnarray}\label{multiplicacao}
B \odot C & = &  \Phi(\Phi^{-1}(B)\cdot\Phi^{-1}(C)) \\\nonumber
 & = &  (rs - qp) + (sq+rp)A, 
\end{eqnarray}
and
\begin{eqnarray}\label{eq:division}
B \odiv C & = &  \Phi(\Phi^{-1}(B)\div\Phi^{-1}(C)) \\\nonumber
& = &   \Phi\left(\frac{rs + qp}{s^2+p^2} + \frac{qs-rp}{s^2+p^2}i\right)\\\nonumber
 & = &  \frac{rs + qp}{s^2+p^2} + \frac{qs-rp}{s^2+p^2}A, 
\end{eqnarray}
provided that $C\neq 0$.  For notational convenience, we also denote the quotient $B \odiv C$ by the usual symbol $\frac{B}{C}$.

In analogy with the so-called complexification process (see \cite{arnold1992ordinary,piccione2008student}) and in view of Theorem \ref{col_ev} and of the fact that $\R \subset \RFA$, the vector space $\RFA$ equipped with the product given by  \eqref{multiplicacao} could be called ``fuzzyfication'' of the vector space $\R$. % (or fuzzy-space $\RFA$ for short).

\end{remark}

It is worth noting that the $\odot$ multiplication operator extends the scalar product. More precisely, let $r \in \R \subset \RFA$ and $B = (s + pA) \in \RFA$, from \eqref{multiplicacao}, we have that 
$$
r \odot B = (r +0A) \odot (s + pA) =  rs + rpA = r(s+pA) = rB.
$$
%Moreover, the powers of $A$ can be defined as follows: 
%$A^0 = 1$, $A^1 = A$, and $A^j = A \odot A^{j-1}$ for all integer $j \geq 2$. 
%From Equation \eqref{multiplicacao}, we also obtain the following expression for the powers of $A$:  
%$A^0 = 1$, $A^1 = A$, $A^2 = -1$, $A^3 = -A$ and $A^j = A^{j \;\mathrm{mod}\; 4}$ for all $j \in [4, \infty) \cap \mathbb{N}$. Hence, 
%the powers of $A$ cycle through $A$, $-1$, $-A$, $1$ because $A^2= -1$. This is analogous to what happens with the powers of $i$ in $\C$.

Since every structure in $\C$ can be carried to the space $\RFA$ through the isomorphism $\Phi$, we can establish a parallel between these two spaces. 
In fact, in the same way that $\C$ is a 2-dimension vector space over real field ($\R$) with $\{1,i\}$ as the canonical basis, the vector space $\RFA$ over $\R$ is $2$-dimension with $\{1,A\}$ as a basis. 
Recall that every field is a 1-dimension over itself. 
Thus, the set of complex numbers is a $1$-dimension vector space over itself, that is, over the complex field ($\C$), and its canonical basis is $\{1\}$. Similarly, using the isomorphism between $\C$ and $\RFA$, we also obtain that 
$\RFA$ is a field and, hence, it is an $1$-dimension vector space over itself (i.e. the field $\RFA$) and, in this case, $\{1\}$ is a basis for $\RFA$. 
This last result is stated in the following theorem.

\begin{theorem}\label{thm:rfa_field}
The space $(\RFA,\oplus,\ominus,\odot,\odiv)$ is field isomorphic to $\C$.
\end{theorem}
\begin{proof}
It follows from the fact that all these operations on $\RFA$ are induced by the isomorphism $\Phi$.  
Note that, from Equation \eqref{eq:division}, for every $C=(s + pA) \neq 0$, we have that 
\[
\frac{1}{C} = \frac{s}{s^2+p^2} + \frac{-p}{s^2+p^2}A.
\]

\end{proof}
It is worth noting that $B \circledast C=\Phi(\Phi^{-1}(B) \ast \Phi^{-1}(C))$, where $\ast$ is one of the arithmetic operations on the field $\mathbb{C}$.

\subsection{Interpretations and Representations}

 In the context of uncertainty modeling, the set $\RFA$ is composed of uncertain quantities given by fuzzy numbers that are linearly correlated with each other since each element is linearly correlated with $A$ \cite{de2021differential}. Additionally, we can treat each element as fuzzy noise around the expected value ``$r$." In statistical terms, ``$r+qA$" is a typical additive uncertainty model where ``$r$" is a population expected value, and ``$qA$" is an individual variability term. In this sense, according to \cite{couso2014statistical}, we can say that we are facing an epistemic approach. This approach aligns with the interpretation provided by Dubois. 
 
 On the other hand, due to the isomorphism between the set of complex numbers and the $\RFA$ space, with $A$ being asymmetric, all interpretations and concepts of complex numbers can be transferred to the space $(\RFA, \oplus,\odot, \|\cdot\|_{\Phi})$. As in $\C$ (and here too), many concepts are purely abstract, such as the derivative at a point. It does not have a geometric meaning like a slope.

\subsubsection{Polar coordinates}

\begin{figure}
\centering
%\begin{overpic}[width=0.85\textwidth,grid,tics=10]{vector_rfa.png}
\begin{overpic}[width=0.85\textwidth,tics=10]{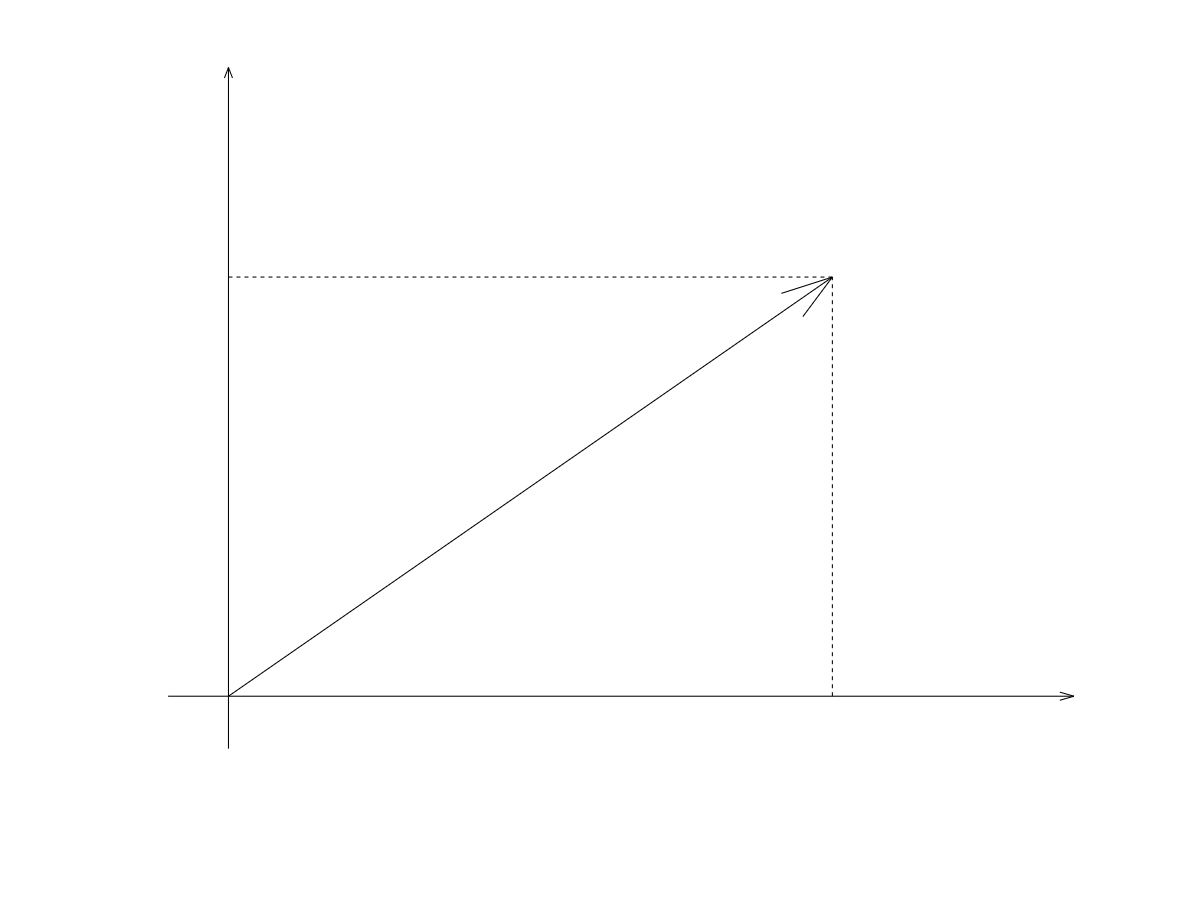}
 \put (70,52) {$\displaystyle z = r +qA$}
 \put (69,14) {$\displaystyle r$}
 \put (16,52) {$\displaystyle q$}
 \put (41,36){\rotatebox{43}{$\displaystyle \| z \|_{\Phi}$}}
 \put (90,13) {$\displaystyle Re$}
 \put (12,69) {$\displaystyle Fu$}
 \put (29,17) {\tikz\draw (0,0) arc (0:44:1);}
 \put (33,20) {$\displaystyle \theta$}
\end{overpic}
\caption{Vector and Polar representation of a fuzzy number $B=(r+qA) \in \RFA$.}
\label{fig:polar}
\end{figure}
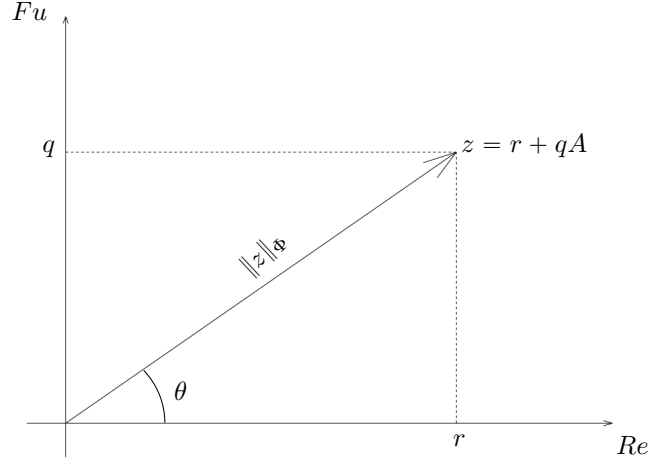

Taking $\RFA$ as a vector space over $\R$ with base $\{1, A\}$ in analogy to the complex case, we can represent $z=r+qA\in\RFA$ as in Figure \ref{fig:polar}.
We denote $r$ by $\text{Re}(z)$ and $q$ by $\text{Fu}(z)$.

The angle $\theta$ between $z$ and $\text{Re}(z)$ is called the argument of $z$ and denoted by $\text{arg }z$. If $r=\|z\|_\Phi\cos{\theta}$ and $q=\|z\|_\Phi\sin{\theta}$, then $z=\|z\|_\Phi(\cos{\theta}+\sin{\theta} A )$. In this case, we have the polar representation of $z$ with $(\|z\|_\Phi, \theta)$ being its polar coordinates.

The conjugate of $z$ is defined by $\overline{z}=r-qA$. One can easily show that  $$z\odot\overline{z}=\|z\|_{\Phi}^2 =r^2+q^2 \quad \mbox{ and } \quad \|z\|_{\Phi}=\sqrt{z\odot\overline{z}}.
$$
Figure \ref{fig:polar} exhibits a Polar representation of an element of $B=r+qA$ of $\RFA$.

Multiplication in polar coordinates has an interesting representation. If $z_1=\|z_1\|_\Phi(\cos{\theta_1}+\sin{\theta_1} A )$ and $z_2=\|z_2\|_\Phi(\cos{\theta_2}+\sin{\theta_2} A )$, then, by \eqref{multiplicacao}, $$z_1\odot z_2=\|z_1\|_\Phi\|z_2\|_\Phi[\cos{(\theta_1+\theta_2)}+\sin{(\theta_1+\theta_2)}A],$$ where
$$\|z_1\odot z_2\|_{\Phi}=\|z_1\|_\Phi\|z_2\|_\Phi.$$ and
$$\\arg(z_1\odot z_2)=arg(z_1)+arg(z_2).$$
We can interpret that in the base $\{1,A\}$, multiplying a fuzzy number $z_1$ by $z_2$ corresponds to performing a rotation and dilation/contraction of $z_1$ according to $\text{arg}(z_2)$ and $\|z_2\|_{\Phi}$, respectively.

Just as with complex numbers, we can rewrite $z=\|z\|_{\Phi}(\cos\theta+\sin\theta A)$ in an exponential form by means of the Euler-type formula: $z=\|z\|_{\Phi}e^{A\theta}$. In addition, we can perform arithmetic operations among the elements of $\RFA$ in polar coordinates using Moivre formula, i.e., $z^n=[r(\cos\theta+\sin\theta A)]^n=r^n(\cos{n\theta}+\sin{n\theta} A)$, while $\sqrt[n]{z}=\sqrt[n]{r}(\cos{\frac{\theta}{n}}+\sin{\frac{\theta}{n}}A)$.    
 To illustrate, consider the fuzzy number $z=1+\sqrt{3}A$ having polar coordinates $\left(2,\frac{\pi}{3}\right)$, thus, we have 
$$\sqrt{z}=\sqrt{2}(\cos{\frac{\pi}{6}}+\sin{\frac{\pi}{6}}A)= \sqrt{2}\left(\frac{\sqrt{3}}{2}+\frac{1}{2} A\right)=B.$$ 
Now, as we mentioned before, real and fuzzy parts are given respectively 
 by $\text{Re}(B) = \sqrt{2}\frac{\sqrt{3}}{2}$ and $\text{Fu}(B) = \sqrt{2}\frac{1}{2}$, however, differently from the complex numbers, they can be joined by the traditional form of a fuzzy number. For example, if $A=(-2;0;4)$, then $B=\sqrt{2}\left(-1+\frac{\sqrt{3}}{2}; \frac{\sqrt{3}}{2}; 2+\frac{\sqrt{3}}{2}\right)$. 

From the isomorphism $\Phi$, we have that the fuzzy number $A$ in $\RFA$ plays a similar role to the imaginary unit ($i$) in $\C$. But, unlike the imaginary components of complex numbers that are not represented in the set of real numbers, the uncertain components, which are associated with the fuzzy number $A$,  
%that plays the role of the imaginary unit by the bijection $\Phi$, 
are represented by fuzzy numbers on $\R$. Thus, an element $(r + qA) \in \RFA$ can be interpreted as a ``fuzzyfication'' of the deterministic/predictable value $r$. 
But, unlike complex numbers, here there is no need to add an additional axis to graphically represent a component of an element of $\RFA$.
This is the main difference in $\C$ and $\RFA$. 

%\begin{comment}
%\textcolor{blue}
%\begin{definition}{(Euler's formula)} Motivated by the expression in the last equality bellow, by Euler's formula (in $\C$) and the fact that $\Phi$ is an isomorphism, we define $e^z$, for any $z=r+qA$ in $\RFA,$ by the formula
%\begin{eqnarray*}
%e^z = e^{r+qA}=e^r(\cos(q) + \sin(q)A),
%\end{eqnarray*}
%and the above formula is called Euler's formula in $\RFA.$\
%\end{definition}

%{\bf Remark:} We emphasize that the value of $e^{z},$ given here by Euler's formula does not coincide with that given by Zadeh extension. By the way, in this last case, $e^{A}$ do not belong to $\RFA,$ if $A$ be triangular fuzzy number. Also, the values for power ($z^n$) as well as $sin(z)$, $cos(z)$ among others, do not coincide with those obtained via Zadeh Extension. This fact is due to the arithmetic operations that make $\RFA$ a Banach space. These arithmetic operations, obtained via $\Phi$ isomorphism, are interactive, while the Zadeh Extension principle uses non-interactive fuzzy arithmetic \cite{carlsson2004additions}.

%\begin{example}
%For $z= 4 + 3A$, by Definition \ref{def:Exponetial}, the exponential of $z$ is given by 
%\[
%e^z = e^{4 + 3A} = e^4(\cos(3) + \sin(3)A) = e^4\cos(3) + e^4\sin(3)A.
%\]
%\end{example}

%As a consequence of definition above, it easy to see that properties of exponential are true:

%\end{comment}

\subsection{Fuzzy Mappings in $\RFA$} 

A fuzzy mapping $f$ in $\RFA$  is nothing more than a function from $D\subseteq\RFA$ to $\RFA$. 
Let $z=(x+y)A \in \RFA$, with $x,y\in\R$.   
The bijection $\Phi$ ensures that %for every $z\in\RFA$ there exist unique $x,y\in\R$ such that $z=x+yA$ and that 
there exist unique functions $u,v:\bar{D}\to\R$ such that
\begin{equation}
   u(x,y)= \text{Re}(f(x+y A))\quad\text{ and }\quad v(x,y)=\text{Fu}(f(x+yA)),
\end{equation}
for all $(x,y) \in \bar{D} = \Phi^{-1}(D) \subseteq \R^2$. 

Using the bijection $\Phi$, many of the important functions in $\C$ \cite{brown2009complex} can be defined similarly in $\RFA$. More precisely, a mapping $f:\C\to\C$ can be carried to the space $\RFA$ by taking $g = \Phi \circ f \circ \Phi^{-1}: \RFA \to \RFA$. 
%Some fundamental mappings from $\RFA$ to $\RFA$ are as follows:
%{\color{blue} 
For example, if $f:\C\to \C$ is given by $f(r + iq) = e^{r + iq} = e^r\cos(q) + ie^r\sin(q)$ (note that this last equality follows from the Euler's formula in $\C$), 
then $g = \Phi \circ f \circ \Phi^{-1}: \RFA \to \RFA$ is given by 
\begin{eqnarray}%\label{eq:Exponential}
\nonumber g(r + qA) & = & \Phi\left(f(\Phi^{-1}(r + qA))\right) \\ 
\nonumber & = & \Phi\left(f(r + iq)\right) \\ 
\nonumber & = & \Phi\left(e^r\cos(q) + ie^r\sin(q)\right) \\
\nonumber & = & e^r\cos(q) + e^r\sin(q)A \\
\nonumber & = & e^r(\cos(q) + \sin(q)A).
\end{eqnarray}
Based on the expression in the last equality, Euler's formula, and the fact that $\Phi$ is an isomorphism, we define $e^z$, for any $z=r+qA$ in $\RFA$, as follows. 
%Motivated by the expression in the last equality, by Euler's formula and by the fact that $\Phi$ is an isomorphism, we define $e^z$, for any $z=r+qA$ in $\RFA$, as follows. 

\begin{definition}\label{def:Exponetial}
For every $z = (r + qA) \in \RFA$, the exponential of $z$ is defined by the fuzzy number $e^{z} \in \RFA$ given by 
\begin{equation}\label{eq:Exponential}
e^z = e^r(\cos(q) + \sin(q)A).
\end{equation}
Equation \eqref{eq:Exponential} is called Euler's formula in $\RFA$ and the mapping $z \mapsto e^z$ is called 
the {\bf exponential function} in $\RFA$.
\end{definition}
%In the following, we propose some properties of exponential functions, as follows:
Some properties of exponential function follow directly from Euler's formula in $\RFA:$
\begin{itemize}
    \item [i.] $e^{z_1}e^{z_2}=e^{z_1+z_2},\forall z_1,z_2\in\RFA,$
    \item [ii.] $(e^z)^n=e^{nz},\forall z\in\RFA,n\in\mathbb{Z}.$
\end{itemize}
%\textcolor{red}{\begin{example}
%Let $z= 4 + 3A$ with $A=(-1,0,1)$. Then, by Definition \ref{def:Exponetial}, the exponential function is given by 
%\begin{align*}
%e^z &= e^{4 + 3A} = e^4(\cos(3) + \sin(3)A) %\\
%%&= e^4\cos(3) + e^4\sin(3)A\\
%&=\left(e^4\cos(3) - e^4\sin(3),e^4\cos(3),e^4\cos(3) + e^4\sin(3)\right).
%\end{align*}
%\end{example}}
\begin{example}
For $z= 4 + 3A$, by Definition \ref{def:Exponetial}, the exponential of $z$ is given by 
\[
e^z = e^{4 + 3A} = e^4(\cos(3) + \sin(3)A) = e^4\cos(3) + e^4\sin(3)A.
\]
\end{example}

In a similar way to the exponential function, other fundamental mappings from $\C$ to $\C$ can also be defined in $\RFA$ as follows:

%{\bf Exponential function:} For every $z = x+yA \in \RFA$, the exponential function is defined as $\exp(z)=e^{x+yA}=e^{x}(\cos{y}+A\sin{y}).$ Note that it is in agreement with an Euler-type formula.

{\bf Logarithmic function:}
For every $z = x+yA \in \RFA$, the logarithmic function is given by $\ln(z)=\ln{\|z\|_{\Phi}}+A(\theta+2n\pi)$, where $n=0,\pm1,\pm2,\ldots$, $z\neq 0$, ${\|z\|_{\Phi}}$ is the modulus of $z$ and $\theta= \text{arg}(z)$ is the argument (angle) of $z$ with the $x$ axis.

{\bf Power function:} 
Let $a\in\R$. For every $z \in \RFA$, $z\neq 0$, we define $z^a=\exp{(a\ln(z))}$. 

{\bf Polynomial function:}
Given $a_0,\ldots,a_n\in\RFA$. For every $z \in \RFA$, we define the polynomial function of $n$ degree as $f(z)=a_0\oplus a_1\odot z\oplus a_2\odot z^2\oplus\ldots\oplus a_n\odot z^n$.

\begin{example}Consider $z=x+yA\in\RFA$.
    The real and fuzzy part of the function $f(z)=exp(z^2\oplus z)$ are given by $\text{Re}(f(z))=u(x,y)= exp(x^2-y^2+x)\cos{(2xy+y)}$ and $\text{Fu}(f(z))=v(x,y)= exp(x^2-y^2+x)\sin{(2xy+y)}$.
\end{example}

%{\color{blue}
\begin{remark}
For a given $z=r+qA \in \RFA$, the symbol $e^{z}$ refers to the fuzzy number in $\RFA$ as defined by Eq. \eqref{eq:Exponential}, which does not coincide with the Zadeh extension of the exponential function $e^x$ at $z$.  In general, the Zadeh extension of the exponential function $e^x$ applied to a non-crisp element of $\RFA$ results in a fuzzy number that does not belong to $\RFA$. 
This observation also applies to other functions, such as logarithmic, power, polynomial, etc, which do not coincide with their Zadeh extension. This outcome is expected because $\RFA$ is a Banach space with interactive arithmetic operations, as induced by the isomorphism $\Phi$ \cite{de2021differential}.
In contrast, the Zadeh extension principle for arithmetic operations leads to the usual or the non-interactive fuzzy arithmetic  \cite{carlsson2004additions}.
\end{remark}

%\begin{remark}
%For some $z=r+qA \in \RFA$, the symbol $e^{z}$ corresponds to the fuzzy number in $\RFA$ given by Eq. \eqref{eq:Exponential} which does not coincide with the Zadeh extension of the exponential function $e^x$ at $z$.  By the way, in general, the Zadeh extension of the exponential function $e^x$  at a non-crisp element of $\RFA$ is a fuzzy number that does not belong to $\RFA$. 
%Such a comment is also applied to functions such as logarithmic, power, polynomial, etc. that do not coincide with their Zadeh extension. One can expect this fact since  $\RFA$ is a Banach space with the arithmetic operations induced by the isomorphism $\Phi$ for which are interactive  \cite{de2021differential}. 
%In contrast, the Zadeh Extension principle of arithmetic operations leads to the usual or the non-interactive fuzzy arithmetic  \cite{carlsson2004additions}.
%\end{remark}

%}

\section{Calculus for Fuzzy Mapping in $\RFA$}\label{Analyticalfuzzyfunctions}

Once again, such concepts can be ``transferred" from the space of complex numbers to the $\RFA$ space via the $\Phi$ operator. Note that $\RFA$ is now equipped with product and quotient operations between its elements, so it makes sense to discuss the limit,  continuity, differentiability, and integrability of functions in $\RFA$. 

\subsection{Limit and continuity}

\begin{definition}
Let $f:D\subseteq\RFA\to\RFA$ and $z_0\in\RFA$. We say that there is the limit of $f$ in $z_0$ if there is $L\in\RFA$ such that for each $\epsilon > 0$ there exists $\delta>0$ such that
\begin{equation}
    z\in D, 0<\|z\ominus z_0\|_{\Phi}<\delta \Rightarrow \|f(z)\ominus L\|_{\Phi}<\epsilon.
\end{equation}
In this case, we denote it by 
\begin{equation}
    \lim_{z\to z_0}{f(z)}=L. 
\end{equation}
\end{definition}

Since $\RFA$ is isomorphic to $\C$, whenever it exists, the limit is unique. 
Moreover, it is worth noting that the limit is independent of the way that $z$ approaches $z_0$ in the fuzzy plane ($\RFA$).

\begin{proposition}
Let $f,g:D\subseteq\RFA\to\RFA$ be functions such that $\displaystyle\lim_{z\to z_0}{f(z)}$ and $\displaystyle\lim_{z\to z_0}{g(z)}$ exist. For every $a \in \RFA$, we have 
\begin{itemize}
    \item[i)] $\displaystyle\lim_{z\to z_0}{[a\odot f(z)\oplus g(z)]}=\left(a\odot \lim_{z\to z_0}{f(z)}\right) \oplus \lim_{z\to z_0}{g(z)}$;
     \item[ii)] $\displaystyle\lim_{z\to z_0}{[f(z)\odot g(z)]}= \lim_{z\to z_0}{f(z)} \odot  \lim_{z\to z_0}{g(z)}$;
     \item[iii)] $\displaystyle\lim_{z\to z_0}{[f(z)\odiv g(z)]}= \lim_{z\to z_0}{f(z)}  \odiv  \lim_{z\to z_0}{g(z)}$ if  $\displaystyle\lim_{z\to z_0}{g(z)}\neq 0$.
\end{itemize}    
\end{proposition}
\begin{proof} 
%The proof of this and the next results is provided in the corresponding result (refer to \cite{brown2009complex}) and the isomorphism $\Phi$.

From isomorphism, $\Phi$, the proof for the $\RFA$ space remains analogous to the case in the complex space $\C$ \cite{brown2009complex}.
\end{proof}

\begin{definition}
 Let $f:D\subseteq\RFA\to\RFA$ and $z_0\in\RFA$. We say that $f$ is continuous in $z_0$ if
\begin{equation}
    \lim_{z\to z_0}{f(z)}=f(z_0). 
\end{equation}
We say that $f$ is continuous in $D$ if $f$ is continuous at all points of $D$.
\end{definition}

Note that the notion continuity can alternatively be defined in terms of open sets since $\RFA$ is also a topology space induced by the norm $\| \cdot \|_\Phi$. Recall that a subset of $Y$ of $\RFA$ is said to be open if for every $B \in Y$ there exists $\epsilon > 0$ such that for all $C \in \RFA$ satisfying $\| B - C \|_\Phi < \epsilon$ implies that $C \in Y$ \cite{rudin64}.

\begin{proposition}\label{limitedaspartesrealefuzzy}
Let $f : D \subseteq \RFA\to\RFA$, $u(x, y) =\text{Re}(f(z))$ and $v(x, y) = \text{Fu}(f(z))$ and $z_0 = x_0+y_0 A \in \RFA$, $x_0, y_0\in\R$. For the limit of $f$ to exist at $z_0$, it is necessary and sufficient that there are limits
of $u$ and $v$ in $(x_0, y_0)$. Then, we have
\begin{equation}
  \lim_{z\to z_0}{f(z)} = \lim_{(x,y)\to(x_0,y_0)}{u(x, y)} +   \lim_{(x,y)\to(x_0,y_0)}{v(x, y)}\, A.
\end{equation}
\end{proposition}
\begin{proof}Immediate consequence of this property in $\C$ and the isomorphism $\Phi$.\end{proof}

\subsection{Derivatives}

\begin{definition}\label{def:derivative}
Let $D \subseteq \RFA$ be an open set, $f : D \to \RFA$ and $z_0\in D$. The derivative of $f$ at $z_0$ is the limit
%We say that $f$ is derivable in $z_0$ if the following limit exists
\begin{equation}
f'(z_0) =\lim_{z\to z_0}{\frac{f(z) \ominus f(z_0)}{z \ominus z_0  }}
\end{equation}
when it exists. Moreover, the derivative $f'(z_0)$ is also denoted by $\frac{df}{dz}(z_0)$
\end{definition}

The basic formulas used to find the derivatives below can be obtained through the complex case \cite{brown2009complex} and the isomorphism $\Phi$.
\begin{proposition}\label{regrasderivada}
    Let $D \subseteq\RFA$ be open, $z_0 = x_0 + y_0 A \in D$, $f, g : D\to\RFA$, and $a \in\RFA$. If $f$ and $g$ are differentiable in $z_0$ then
    \begin{itemize}
        \item[i)]$\displaystyle (a\odot f\oplus g)'(z_0)=a\odot f'(z_0)\oplus g'(z_0)$;
        \item[ii)]$\displaystyle (f\odot g)'(z_0)=f(z_0)\odot g'(z_0)\oplus f'(z_0)\odot g(z_0)$;
      \item[iii)]$\displaystyle (f\odiv g)'(z_0)=  \left(f'(z_0)\odot g(z_0)\oplus f(z_0)\odot g'(z_0)\right) \odiv g(z_0)^2$, when $g(z_0)\neq 0$.  
    \end{itemize}
\end{proposition}
\begin{proof}
It follows from the properties of complex calculus and from the isomorphism $\Phi$.
\end{proof}

\begin{example} \label{ex:derivatives}
Consider $f:\RFA\to \RFA$ and $z_0\in\RFA$. 
    \begin{enumerate}
        \item If $f(z) = k$, $k\in\RFA$, then $f'(z_0) = 0$.
        \item If $f(z) = z^n$, $n\in\mathbb{N}$, then $f'(z_0) = nz_0^{n-1}$.
        \item If $f(z) = z^{-n}$, $n\in\mathbb{N}$, then $f'(z_0) = -nz_0^{-n-1}$, $z_0\neq 0$.
        \item If $f(z) = e^{z}$, then $f'(z_0) = e^{z_0}$.
        \item If $f(z)=a_0\oplus a_1\odot z\oplus a_2\odot z^2\oplus\ldots\oplus a_n\odot z^n$, then $$f'(z_0)= a_1\oplus 2a_2\odot z_0\oplus\ldots \oplus n a_n\odot z_0^{n-1}.$$
    \end{enumerate}
\end{example}

The following results are consequences of the isomorphism $\Phi$ and, therefore, the proofs are analogous to the complex case (see \cite{brown2009complex}).
\begin{proposition}\label{prop:CauchyRiemman}
 Let $D \subseteq\RFA$ be open, $z_0 = x_0 + y_0 A \in D$, $x_0, y_0\in\R$, $f : D \to \RFA$, $u=u(x, y) =\text{Re}f(z)$ and $v=v(x, y) = \text{Fu}f(z)$. If $f'(z_0)$ exists, then the partial derivatives of $u$ and $v$ in ($x_0, y_0$) exist and satisfy the Cauchy-Riemann equations:
\begin{equation}\label{condicoesCR}
\frac{\partial u}{\partial x}(x_0, y_0)=\frac{\partial v}{\partial y}(x_0, y_0)\quad\text{ and }\quad \frac{\partial u}{\partial y}(x_0, y_0)=-\frac{\partial v}{\partial x}(x_0, y_0).
\end{equation}
Furthermore, 
\begin{equation}\label{CRsolution}
    f'(z_0)=\frac{\partial u}{\partial x}(x_0, y_0)+ \frac{\partial v}{\partial x}(x_0, y_0)A=\frac{\partial v}{\partial y}(x_0, y_0)- \frac{\partial u}{\partial y}(x_0, y_0)A.
\end{equation}
\end{proposition}
\begin{theorem}\label{teoCR}
Let $D \subseteq\RFA$ be open, $z_0 = x_0 + y_0 A \in D$, $x_0, y_0\in\R$, and $f : D \to \RFA$ a
function such that, $u=u(x, y) =\text{Re}f(z)$ and $v=v(x, y) = \text{Fu}f(z)$ have continuous first-order partial derivatives in $(x_0, y_0)$. If $u$ and $v$ satisfy the Cauchy-Riemann equations \eqref{condicoesCR} then there exists $f'(z_0)$ and it is given by \eqref{CRsolution}.
\end{theorem}

\begin{example}
  Consider $z=x+yA\in\RFA$ and $f:\RFA\to\RFA$, such that %$$f(z)=2A\odot_A e^{\lambda (A\ominus_A z)}=2A\odot_A\exp{[\lambda (-x+(1-y)A)]}.$$ 
    \begin{eqnarray}\label{ex1}
    f(z)=e^{- z}= e^{- x}\left(\cos{y}- \sin{y}A\right).
    \end{eqnarray}
    We can rewrite $\eqref{ex1}$ separating the real and fuzzy parts as follows
    \begin{eqnarray*}
        u(x,y)=e^{- x}\cos{ y}\qquad \text{ and } \qquad v(x,y)=-e^{- x}\sin{ y}
    \end{eqnarray*}
    since the above functions satisfy \eqref{condicoesCR} everywhere and since these derivatives are everywhere
continuous, the conditions in the Theorem\ref{teoCR} are satisfied. Thus $f'(z)$ exists everywhere, and 
    \begin{equation}
        f'(z)= \left(- e^{- x}\cos{y}\right) + \left(- e^{- x}\sin{y}\right)A\textcolor{blue}{.}
    \end{equation}
    
    \end{example}

\begin{remark}
Calculating this derivative does not involve the analysis of $\alpha$-levels. However, we could write $$f(z)=f(x+yA)=u(x,y)+v(x,y)A$$ as $$[f(z)]_{\alpha}=[f(x+Ay)]_{\alpha}=u(x,y)+ v(x,y)[A]_{\alpha},$$ 
for all $\alpha \in [0,1]$. In this way, by Cauchy-Riemann conditions, we have $$[f'(z)]_{\alpha}=[f'(x+Ay)]_{\alpha}=\frac{\partial u}{\partial x}(x,y)+ \frac{\partial v}{\partial x}(x,y)[A]_{\alpha},$$
or equivalently
$$[f'(z)]_{\alpha}=[f'(x+yA)]_{\alpha}=\frac{\partial v}{\partial y}(x,y)+ \left(-\frac{\partial u}{\partial y}(x,y)\right )[A]_{\alpha},$$
for all $\alpha \in [0,1]$. Thus,  the change in the sign of $\frac{\partial v}{\partial x}(x,y)$ (or equivalently $-\frac{\partial u}{\partial y}(x,y)$) is what determines the appearance of switching points \cite{bede13}. 
\end{remark}
The next definition is adapted from the complex case \cite{brown2009complex}.
\begin{definition}
    Let $D$ be open, $z_0\in D\subseteq\RFA$ and $f : D \to\RFA$. We say that $f$ is analytic in $z_0$ if $f$ is derivable at all points of some open disk centered in $z_0$. We say that $f$ is analytic in $D$ if $f$ is analytic at all points of $D$. %An analytic function in C is called an integer function.
\end{definition}

\begin{proposition}[Chain rule]\label{prop:chain_rfa}
    Let $D,\Omega \subseteq \RFA$  be open sets, $f : D \to \Omega$ and $g :\Omega\to\RFA$. If $f$ is analytic in $\Omega$ and $g$ is analytic in $\Omega$ then the composite $g \circ f : D\to\RFA$ is also analytic in $D$ and holds 
    \begin{equation}\label{eq:complex_chain_rule}
    (g \circ f)'(z_0) = g'(f(z_0))\odot f'(z_0), \forall z_0\in D.    
    \end{equation}    
\end{proposition}
\begin{proof}
It follows from the bijection $\Phi$ and from the well-known results from the complex calculus theory.   
\end{proof}

As an immediate consequence of Proposition \ref{prop:chain_rfa} and Example \ref{ex:derivatives}, we obtain the following example. 

\begin{example} \label{ex:derivative_exp}
The function $f:\RFA\to \RFA$ given by $f(z) = e^{k\odot z}$, with $k\in\RFA$, is differentiable and $f'(z_0) = k\odot e^{k\odot z_0}$ for every $z_0 \in \RFA$.
\end{example}

\subsection{Antiderivative Function}%Integrability
 The following results are adapted from the complex case \cite{brown2009complex}.
\begin{definition}
    Let $f : \Omega\subseteq\RFA\to\RFA$. We say that $F : \Omega \to\RFA$ is an antiderivative (primitive) of $f$ if $F'(z) = f(z)$ for all $z\in\Omega$.
\end{definition}
In this case, it is common to denote $F$ by $\int{f(z)dz}$.
\begin{theorem}\label{teoIntegral}
    If $\Omega\subseteq\RFA$ is a simply connected set and $f$ is an analytical function then, fixed $z_0\in \Omega$ there exists $C\in\RFA$ such that the general form of a primitive $f$ is given by
    \begin{equation}\label{integral1}
        F(z)=\int_{z_0}^{z}{f(\zeta)d\zeta+C},\quad z\in\Omega.
    \end{equation}
    %is a primitive of $f$.
    Given $z_0$, formula \eqref{integral1} indicates that the integral of $f$ is independent of the path that connects $z_0$ to $z$.
\end{theorem}
%\begin{proof}
  %  Via isomorphism $\Phi$, the proof of the theorem remains analogous to the case of complexes \cite{brown2009complex}.
%\end{proof}
From Theorem \ref{teoIntegral}, it immediately follows that 
\begin{equation}
    \int_{z_1}^{z_2}{f(z)dz}=F(z_2)\ominus F(z_1).
\end{equation}

\begin{example}
Consider $f:\RFA\to\RFA $. The continuous function $f (z) = z^2$ has an antiderivative $F(z) = \displaystyle\frac{z^3}{3}$ throughout the plane. Thus 
\begin{equation}
    \int_{0}^{1+A}{z^2 \,dz}=\frac{z^3}{3}\bigg{]}_{0}^{1+A}=\frac{1}{3}(1+A)^3= -\frac{2}{3}+ \frac{2}{3}A.
\end{equation}
\end{example}

\section{Fuzzy curves and ordinary differential equations in $\RFA$} \label{Fuzzyfunctions}

A fuzzy curve in $\RFA$ is a fuzzy-number-valued mapping $w:I\subseteq\R\to\RFA$ given by
\begin{equation}\label{funct}
    w(t)=x(t)+y(t)A,\quad t\in I\subseteq\R.
\end{equation}

Since addition and multiplication are defined in $\RFA$, we can define the pointwise sum and product of fuzzy curves as follows: 
\begin{eqnarray}
    (w_1\oplus w_2)(t)&=&w_1(t)\oplus w_2(t)\quad \text{and}\\\nonumber
    (w_1\odot w_2)(t)&=&w_1(t)\odot w_2(t)\quad \forall t \in I.
\end{eqnarray}

\subsection{Derivatives and Integrals of the fuzzy curve}

The derivative of the function \eqref{funct} at a point $t_0$ is defined as
\begin{equation}\label{dert1}
  \frac{dw}{dt}\bigg{|}_{t=t_0}=\lim_{h\to 0}{\frac{w(t_0+h)\ominus w(t_0)}{h}}.
\end{equation}
or, equivalently, as 
\begin{equation}\label{dert}
\frac{dw}{dt}\bigg{|}_{t=t_0}=w'(t_0)=x'(t_0)+y'(t_0)A,\quad t_0\in I,
\end{equation}
provided the existence of the derivatives $x'$ and $y'$ at $t_0$ \cite{esmi2018frechet}.

The linearity of the derivative operator (\ref{dert1}) with respect to the addition $\oplus$ and the scalar product in the real field follows immediately from the results of  \cite{esmi2018frechet} combined to the standard bijection between $\C$ and $\R^2$. 

In what follows we state some properties with respect to the multiplication $\odot$ given as in Equation \eqref{multiplicacao}.

\begin{theorem}[Chain rule]\label{thm:chain}
Let $\Omega \subseteq \RFA$ be an open set, $w:(a,b) \to \Omega$, and $f:\Omega \to \RFA$. 
If the fuzzy curve $w$ is differentiable at $t_0 \in (a,b)$   
and $f$ is differentiable at $w(t_0)$, then  
the fuzzy curve $\gamma=f\circ w: (a,b)\to \RFA$ is differentiable at $t_0$ with 
\begin{eqnarray}\label{eq:chain_rule}
\gamma'(t_0) = f'(w(t_0))\odot w'(t_0).    
\end{eqnarray}
\end{theorem}
\begin{proof}
Let $x,y:(a,b)\to \R$ and $u,v:D\to\R$ such that 
$w(t) = x(t) + y(t)A$ and $f(p + qA) = u(p,q) + v(p,q)A$, with $D = \Phi^{-1}(\Omega)$. 
Note that, from the bijectivity of $\Phi$, these functions are unique. 

On the one hand, we have 
\begin{eqnarray*}
\gamma(t) & = & f(w(t)) \\
 & = & f(x(t) + y(t)A) \\
 & = & u(x(t),y(t)) + v(x(t),y(t))A.
\end{eqnarray*}
Using \eqref{dert} and chain rule for real functions, we obtain  
\begin{eqnarray*}
\gamma'(t_0) & = &  \left(\frac{\partial u}{\partial x}(x(t_0),y(t_0))x'(t_0) + \frac{\partial u}{\partial y}(x(t_0),y(t_0))y'(t_0) \right)  + \\
& = &  \left(\frac{\partial v}{\partial x}(x(t_0),y(t_0))x'(t_0) + \frac{\partial v}{\partial y}(x(t_0),y(t_0))y'(t_0) \right)A. \\
\end{eqnarray*}
This last equality proves that $\gamma$ is differentiable at $t_0$ since every derivative and partial derivative and right side exist. 

On the other hand, we have that $w'(t_0) = x'(t_0) + y'(t_0)A$ and, from Proposition \ref{prop:CauchyRiemman}, we have  
$$
f'(x(t_0) + y(t_0)A) = 
\frac{\partial u}{\partial x}(x(t_0),y(t_0))+ \frac{\partial v}{\partial x}(x(t_0),y(t_0))A,
$$
and the Cauchy-Riemman equations is satisfied at $x(t_0) + y(t_0)A$. 

Thus, we have
\begin{eqnarray*}
f'(w(t_0))\odot w'(t_0) & = & \left(\frac{\partial u}{\partial x}(x(t_0),y(t_0))+ \frac{\partial v}{\partial x}(x(t_0),y(t_0))A\right) \odot  \left(x'(t_0) + y'(t_0)A\right) \\
 & = & \left(\frac{\partial u}{\partial x}(x(t_0),y(t_0))x'(t_0) -\frac{\partial v}{\partial x}(x(t_0),y(t_0))y'(t_0)\right) +    \\
 &   & \left(\frac{\partial v}{\partial x}(x(t_0),y(t_0))x'(t_0) + \frac{\partial u}{\partial x}(x(t_0),y(t_0))y'(t_0)\right)A.
\end{eqnarray*}
From the Cauchy-Riemann conditions, it follows that 
\begin{eqnarray*}
f'(w(t_0))\odot w'(t_0) & = &  \left(\frac{\partial u}{\partial x}(x(t_0),y(t_0))x'(t_0) + \frac{\partial u}{\partial y}(x(t_0),y(t_0))y'(t_0) \right)  + \\
&  &  \left(\frac{\partial v}{\partial x}(x(t_0),y(t_0))x'(t_0) + \frac{\partial v}{\partial y}(x(t_0),y(t_0))y'(t_0) \right)A. \\
\end{eqnarray*}
Therefore, we conclude that \eqref{eq:chain_rule} holds true. 

\end{proof}    

\begin{example}
    From Examples \ref{ex:derivatives} and \ref{ex:derivative_exp} and Theorem  \ref{thm:chain}, it follows that
\begin{enumerate}
    \item $\displaystyle\frac{d}{dt}[k\odot w(t)]=k\odot w'(t)$, $\forall k\in\RFA$;
    \item  $\displaystyle\frac{d}{dt}[e^{kt}]=k\odot e^{kt}$, $\forall k\in\RFA$;
\end{enumerate}
\end{example}

The definite integral of $w(t)$ over an interval $a \leq t \leq b$ is defined as
\begin{equation}\label{intt}
     \int_{a}^{b}{w(t) dt}=\int_{a}^{b}{x(t) dt}+\left(\int_{a}^{b}{y(t) dt}\right)A,
\end{equation}
provided that the individual integrals on the right exist.
\begin{example}
    From Equations \eqref{multiplicacao} and  \eqref{intt}, it follows that
\begin{eqnarray}
    \int_{0}^{1}{(1+tA)^2 dt}&=&\left(\int_{0}^{1}{1-t^2 dt}\right)+ \left(\int_{0}^{1}{2t dt}\right) A\\\nonumber
    &=& \frac{2}{3}+ A.
\end{eqnarray}
\end{example}

\subsection{Linear Ordinary Differential Equations}\label{LinearDifferential}

Consider the linear differential equation given by
\begin{equation}\label{lineqfuzz}
\left\{ 
\begin{array}{rcl }
   \frac{d}{dt}w(t)&=&\lambda\odot w(t)\\
   w(0)&=&w_0
   \end{array}
,\right.
\end{equation}
where $\lambda\in\RFA$ and $w:I\subset\R\to\RFA$.

By computing its derivative, one can verify that the curve 
\begin{equation}\label{sollineqfuzz}
w(t)=w_0\odot e^{\lambda t}
\end{equation}
is a solution of \eqref{lineqfuzz}.

In what follows, we present some considerations about the fuzzy function $e^{\lambda t}$. % of real variable $t$. 
Here, in analogy with the classical case, for every fixed $t$, the mapping  $\lambda \to e^{\lambda t}$ is called the phase flow of Equation  \eqref{lineqfuzz}. 

Consider $\lambda=\lambda_1+\lambda_2 A\in\RFA$, where $\lambda_1,\lambda_2\in\R$.
\begin{enumerate}
    \item[a.] If $\lambda_1\neq 0$ and $\lambda_2=0$, then we have that $e^{\lambda t}=e^{\lambda_1 t}\in \R$. In this case, the phase flow of \eqref{lineqfuzz} consists of dilatations/contractions by the factor $e^{\lambda_1 t}$.
    \item[b.]If $\lambda_1=0$ and $\lambda_2\neq 0$, we have that $e^{\lambda t}=e^{(\lambda_2 A) t}\in \RFA$. In this case, $$e^{(\lambda_2 A) t}=\cos{(\lambda_2 t)}+\sin{(\lambda_2 t)}A.$$ Thus, the phase flow of \eqref{lineqfuzz} is a family of oscillations $\{e^{(\lambda_2 A) t}\}$ with angular frequency equal to  $\lambda_2$.
    \item[c.]If $\lambda_1\neq0$ and $\lambda_2\neq 0,$ then we have that $e^{\lambda t}=e^{(\lambda_1+\lambda_2 A) t}\in \RFA$. In this case, $$e^{(\lambda_1+\lambda_2 A) t}=e^{\lambda_1 t}[\cos{(\lambda_2 t)}+\sin{(\lambda_2 t)}A].$$ Thus, the phase flow of \eqref{lineqfuzz} is a dilatation/contraction by a factor $e^{\lambda_1 t}$ with  simultaneous oscillations of angular frequency equal to  $\lambda_2$.
\end{enumerate}

    Alternatively, assuming $w(t) = x(t) + y(t)A$ and $w_0 = x_0 + y_0A$, it is also possible to analyze the fuzzy system \eqref{lineqfuzz} from the following 2-dimensional real system:

\begin{equation}\label{lineqreal}
    \left(\begin{array}{c}
        x^{'}(t) \\
        y^{'}(t)
    \end{array}\right)=
     \left(\begin{array}{cc}
    \lambda_1   & -\lambda_2\\
      \lambda_2 &\lambda_1
   \end{array}\right)
    \left(\begin{array}{c}
        x(t) \\
        y(t)
    \end{array}\right),
\end{equation}
%where $x(t_0)=x_0$ and $y(t_0)=y_0$.
whose solution is given by
\begin{equation}\label{sollineqreal}
\left( 
\begin{array}{c}
 x(t)      \\
 y(t)
\end{array}
\right)=x_0 e^{\lambda_1 t} 
\left( 
\begin{array}{c}
 \cos{(\lambda_2 t)}  \\
\sin{(\lambda_2 t) }
\end{array}
\right)
+y_0 e^{\lambda_1 t}
\left( 
\begin{array}{c}
 -\sin{(\lambda_2 t)}  \\
\cos
{(\lambda_2 t)} 
\end{array}
\right).
\end{equation}
Therefore, the $\text{Re}(w(t))$ and  $\text{Fu}(w(t))$ of \eqref{sollineqfuzz} are given by $x(t)$ and $y(t)$ in \eqref{sollineqreal}. In fact, the variation of the real and fuzzy parts of \eqref{sollineqfuzz}  with time are seen as linear combinations of $e^{\lambda_1 t}\cos{(\lambda_2 t)}$ and $e^{\lambda_1 t}\sin{(\lambda_2 t)}$. 

Given the equations $x(t)$ and $y(t)$ in \eqref{sollineqreal}, the general solution \eqref{sollineqfuzz} of  \eqref{lineqfuzz} can be obtained as follows:
\begin{eqnarray}\label{soljuntas}
    w(t)&=&x(t)+y(t)A\\\nonumber
    &=&[x_0e^{\lambda_1t}\cos{(\lambda_2 t)}-y_0e^{\lambda_1t}\sin{(\lambda_2 t)}]\\\nonumber
    &+&[y_0e^{\lambda_1t}\cos{(\lambda_2 t)}+x_0e^{\lambda_1t}\sin{(\lambda_2 t)}]A\\\nonumber
    &=&(x_0+y_0A)\odot e^{(\lambda_1+\lambda_2A)t}=w_0\odot e^{\lambda t}.
\end{eqnarray}

If $[A]_{1}=\{0\}$,  then $[w(t)]_1=\{x(t)\}$ for all $t$. Thus, the solution can be seen as the ``composition'' of $x(t)$ and the fuzzy part $y(t) = (y_0\cos{(\lambda_2 t)}+x_0\sin{(\lambda_2 t)})A$. In this case, $y$ can be interpreted as the fuzzy noise around the deterministic curve (or the drift) $x(t)$.

Notice that the phase flow of \eqref{lineqreal} with nonreal eigenvalue $\mu=\lambda_1\pm i\lambda_2$ is similar to a family of dilations/contractions by factors $e^{\lambda_1t}$ with simultaneous oscillations with angular frequency of  $\lambda_2$.
As in \cite{arnold1992ordinary}, we refer to the process of transforming \eqref{lineqfuzz} into \eqref{lineqreal} as  ``realification'' of the fuzzy system. We refer to the opposite process as ``fuzzification''.

In the following subsection, we compare the study above with the approach found in \cite{laiate2021bidimensional, laiate2021cross}, where the authors introduced a different notion of multiplication in $\RFA$ named $\Psi$-cross product. To this end, let us analyze the differential equation given in \eqref{lineqfuzz} replacing the multiplication ($\odot$) by the so-called $\Psi$-cross product. This last concept is only defined for fuzzy numbers with 1-$level$ containing only one element (\cite{laiate2021bidimensional, laiate2021cross}).

\subsubsection{Linear Differential Equation in $\RFA$ with $\Psi$-cross product}

First of all, let us recall the definition of the $\Psi$-cross-product \cite{laiate2021bidimensional, laiate2021cross}.

\begin{definition}\label{def:cross}
    Let $A$ be an asymmetric fuzzy number such that $[A]_1=\{a\}$, that is, $[A]_1$ is a unitary set. 
    The $\Psi$-cross product between $B\in\RFA$ and $C\in\RFA$ is given by
\begin{equation}\label{crossproductdef}
	P=B\odot_{\Psi}C=c_1B \oplus b_1C \ominus b_1c_1,
\end{equation}
where $[B]_1=\{b_1\}$ and $[C]_1=\{c_1\}$. 
\end{definition}

In this subsection, unless otherwise stated, we assume that the fuzzy number $A$ satisfies the conditions of Definition \ref{def:cross}. If $B=r_B+q_BA$ and $C=r_C+q_CA$, then 
$b_1 = r_B+q_Ba$, $c_1= r_C+q_Ca$, and 
$$B\odot_{\Psi}C =  (r_Br_C - a^2q_Cq_B) +  [q_B(r_C+q_Ca) + q_C(r_B+q_Ba)]A$$
where $[A]_1 = \{a\}.$ 
For example, if $a=0$, $B=1+2A$, and $C=2+3A$, then $B\odot_{\Psi}C=2B+_{\Psi}1C-1 \times 2=2+7A$.

 Consider the following linear system
\begin{equation}\label{crescimento}
\left\{ 
\begin{array}{c }
w^{'}(t)=\lambda \odot_\Psi w(t) \\
w(0)=w_0\in \mathbb{R}_{\mathcal{F}(A)}
\end{array}
,\right. 
\end{equation}
where $\lambda=\lambda_1+\lambda_2A\in\RFA$,  $w(t)\in\mathbb{R}_{\mathcal{F}(A)}$, and $\odot_{\Psi}$ is the $\Psi$-cross product.  The realification of \eqref{crescimento} is given by
\begin{equation}\label{crescimentoreal}
    \left(\begin{array}{c}
        x^{'}(t) \\
        y^{'}(t)
    \end{array}\right)=
     \left(\begin{array}{cc}
    \lambda_1   & -a^2\lambda_2\\
      \lambda_2 &\lambda_1+2a\lambda_2
   \end{array}\right)
    \left(\begin{array}{c}
        x(t) \\
        y(t)
    \end{array}\right).
\end{equation}
Since both eigenvalues associated to \eqref{crescimentoreal} are equal to $\lambda_1+a\lambda_2$, the solution of \eqref{crescimento} is given by
\begin{equation}\label{solcrescimentoreal}
\left( 
\begin{array}{c}
 x(t)      \\
 y(t)
\end{array}
\right)=y_0 e^{(\lambda_1+a\lambda_2 )t} 
\left( 
\begin{array}{c}
-a  \\
1
\end{array}
\right)
+(x_0+ay_0) e^{(\lambda_1+a\lambda_2 ) t}\left[
t\left( 
\begin{array}{c}
 -a \\
1 
\end{array}
\right)+\left( 
\begin{array}{c}
 1  \\
0
\end{array}
\right)\right]
\end{equation}

Considering that $A$ is an asymmetric fuzzy number centered on zero, that is, $[A]_1=\{a\}=\{0\}$, \eqref{solcrescimentoreal} becomes
\begin{equation}\label{solcrescimentoreal2}
\left( \begin{array}{c} x(t)      \\ y(t)\end{array}\right)= e^{\lambda_1 t} \left( \begin{array}{c}x_0  \\y_0\end{array}\right)+
t e^{\lambda_1 t}
\left( \begin{array}{c} 
0  \\
x_0
\end{array}\right).
\end{equation}
Thus, the fuzzy solution is given by 
\begin{equation}\label{psisol}
w(t) = x(t)+y(t)A = x_0e^{\lambda_1t}+(y_0+x_0t)e^{\lambda_1t}A.
\end{equation}
Note that $w(t)$ does not depend on $\lambda_2$ and, therefore, the uncertainty in the parameter $\lambda$ is completely neglected in the obtained solution. In other words, we get the same solution for $\lambda_2 =0$ (i.e, when $\lambda$ is a real number) or for any other value of $\lambda_2$.  

Some differences between the solution \eqref{psisol} and that given in \eqref{soljuntas} (using the multiplication $\odot$) deserve to be highlighted. In equation \eqref{psisol}, the real part $x(t)$ (the drift) increases or decreases monotonically with time, depending only on the sign of $\lambda_1$. In contrast,  in \eqref{soljuntas}, $x(t)$ oscillates periodically with time, except for $\lambda_2 =0.$ Moreover, according to \eqref{psisol}, the fuzzy part $y(t)$ (or the fuzzy noise) increases or decreases monotonically around zero with time, depending on signs of $\lambda_1$, $x_0$ and $y_0$. On the other hand, according to equation \eqref{soljuntas}, $y(t)$ oscillates periodically with time, with $|y(t)|$ ranging from zero (where the uncertainty collapses) to a certain non-zero maximum value. For example, if $x_0=y_0$, then $|y(t)|=0$ and $w(t)=\pm\sqrt{2}x_0e^{\lambda_1t}$ whenever $\lambda_2t=\frac{3\pi}{4} +k\pi$. On the other hand, when $\lambda_2t=\frac{\pi}{4} +k\pi$, we have $w(t)=\pm\sqrt{2}y_0e^{\lambda_1t}A$.
In all cases, except $\lambda_1 = \lambda_2 = 0$, $w=0$ is the unique equilibrium point which is stable for $\lambda_1<0$, 
since $w(t)\to 0$ as $t\to \infty$ if $\lambda_1<0$, and unstable for $\lambda_1>0.$\

\begin{figure}[H]
\center
\includegraphics[width=0.55\textwidth]{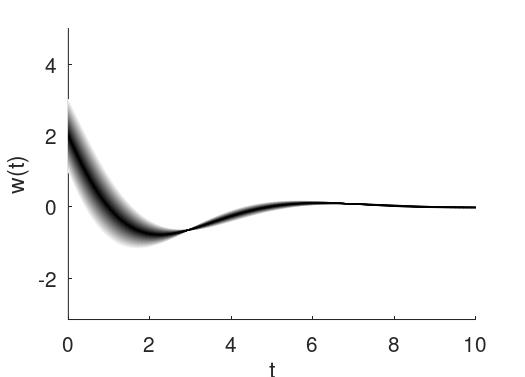}
\caption{Graphical representation of the fuzzy curve $w$ given by Eq. \eqref{soljuntas}, which is the solution of linear fuzzy differential equation \eqref{lineqfuzz} with $A=(-0.5;0;0.51)$, $\lambda_1=-0.5, \lambda_2=0.8$, and $x_0 = y_0=2$. 
Gray lines, from white to black, represent the endpoints of the $\alpha$-level of $w$ for $\alpha$ varying from 0 to 1, respectively.}
\label{fig:decay}
\end{figure}

\begin{figure}[H]
\center
\includegraphics[width=0.55\textwidth]{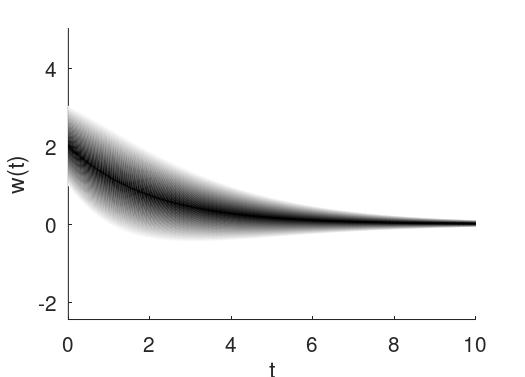}
\caption{Graphical representation of the fuzzy curve $w$ given by Eq. \eqref{psisol}, which is the solution of linear fuzzy differential equation \eqref{crescimento} with $A=(-0.5;0;0.51)$, $\lambda_1=-0.5, \lambda_2=0.8$, and $x_0 = y_0=2$. Gray lines, from white to black, represent the endpoints of the $\alpha$-level of $w$ for $\alpha$ varying from 0 to 1, respectively.}
\label{fig:decay_cross}
\end{figure}

\begin{figure}[H]
\center
\includegraphics[width=0.55\textwidth]{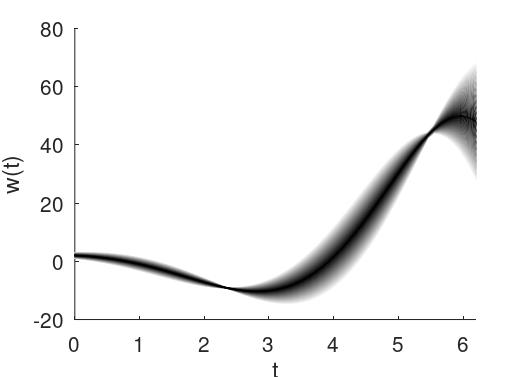}
\caption{Graphical representation of the fuzzy curve $w$ given by Eq. \eqref{soljuntas}, which is the solution of linear fuzzy differential equation \eqref{lineqfuzz} with $A=(-0.5;0;0.51)$, $\lambda_1=0.5, \lambda_2=1$, and $x_0 = y_0=2$. 
Gray lines, from white to black, represent the endpoints of the $\alpha$-level of $w$ for $\alpha$ varying from 0 to 1, respectively.}
\label{fig:grown}
\end{figure}

\begin{figure}[H]
\center
\includegraphics[width=0.55\textwidth]{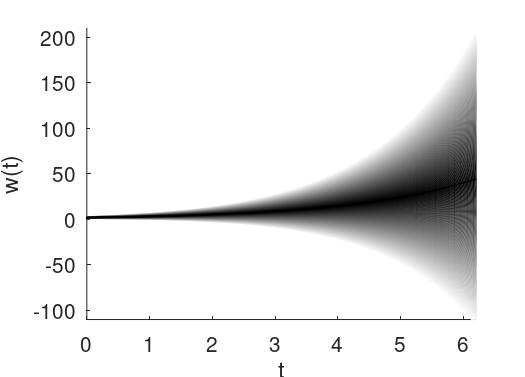}
\caption{Graphical representation of the fuzzy curve $w$ given by Eq. \eqref{psisol}, which is the solution of linear fuzzy differential equation \eqref{crescimento} with $A=(-0.5;0;0.51)$, $\lambda_1=0.5, \lambda_2=1$, and $x_0 = y_0=2$.. Gray lines, from white to black, represent the endpoints of the $\alpha$-level of $w$ for $\alpha$ varying from 0 to 1, respectively.}
\label{fig:grown_cross}
\end{figure}

The comments above can be observed in graphics of Figures \ref{fig:decay} -- \ref{fig:grown_cross}. 
In particluar, Figures \ref{fig:decay} and \ref{fig:decay_cross} exhibit respectively the graphical representation of solutions \eqref{soljuntas} and \eqref{psisol} with $A=(-0.5;0;0.51)$, $\lambda= -0.5 + 0.8(-0.5;0;0.51) = (-0.9;-0.5;-0,092)$, and $w_0= 2 + 2(-0.5;0;0.51) = (1;2;3.02)$. These figures illustrate the case in which $\lambda_1 < 0$ where the solutions $w(t)$ converge to the stable equilibrium point $0$. 
Similarly, Figures \ref{fig:grown} and \ref{fig:grown_cross} exhibit respectively the graphical representation of solutions \eqref{soljuntas} and \eqref{psisol} with $A=(-0.5;0;0.51)$, $\lambda= 0.5 + 1(-0.5;0;0.51) = (0;1;1,01)$, and $w_0= 2 + 2(-0.5;0;0.51) = (1;2;3.02)$. These figures illustrate the case where $\lambda_1 > 0$ which 
phase flows of \eqref{lineqfuzz} and \eqref{crescimento} consist of dilatations by the factor $e^{\lambda_1 t}$. 

%\begin{figure}[H]
%\center
%\includegraphics[width=0.7\textwidth]%{teste_RFA_malthus_prod_complex_level.jpg}
%\caption{Graphical representation of the solution \eqref{soljuntas} of linear fuzzy differential equation \eqref{lineqfuzz}. The dashed lines represent the 0-level and the solid line represents the 1-level. The parameters are $\lambda_1=-0.5, \lambda_2=0.8$, $x_0 = y_0=2$, and the triangular fuzzy number $A=(-0.5;0;0.51)$.}
%\label{fig:decay}
%\end{figure}

%\begin{figure}[H]
%\center
%\includegraphics[width=0.7\textwidth]{teste_RFA_malthus_prod_cross_level.jpg}
%\caption{Graphical representation of the solution \eqref{psisol} of linear fuzzy differential equation \eqref{crescimento}. The dashed lines represent the 0-level and the solid line represents the 1-level. The parameters are $\lambda_1=-0.5, \lambda_2=0.8$, $x_0=y_0=2$, and the triangular fuzzy number $A=(-0.5;0;0.51)$.}
%\label{fig:decay_cross}
%\end{figure}

%\begin{figure}[H]
%\center
%\includegraphics[width=0.7\textwidth]{teste_RFA_malthus_crescimento_prod_complex_level.jpg}
%\caption{Graphical representation of the solution \eqref{soljuntas} of linear fuzzy differential equation \eqref{lineqfuzz}. The dashed lines represent the 0-level and the solid line represents the 1-level. The parameters are $\lambda_1=0.5, \lambda_2=1$, $x_0=y_0=2$, and the triangular fuzzy number $A=(-0.5;0;0.51)$.}
%\label{fig:grown}
%\end{figure}

%\begin{figure}[H]
%\center
%\includegraphics[width=0.7\textwidth]{teste_RFA_malthus_crescimento_prod_cross_level.jpg}
%\caption{Graphical representation of the solution \eqref{psisol} of linear fuzzy differential equation \eqref{crescimento}. The dashed lines represent the 0-level and the solid line represents the 1-level. The parameters are $\lambda_1=0.5, \lambda_2=1$, $x_0=y_0=2$, and the triangular fuzzy number $A=(-0.5;0;0.51)$.}
%\label{fig:grown_cross}
%\end{figure}

We conclude this study of differential equations in $\RFA$ by highlighting that Subsection \ref{LinearDifferential} deals with differential equations whose solution domains are intervals of $\R$, that is, they are curves in $\RFA$.
Due to  the isomorphism $\Phi$,
%Due to the close connection between the spaces $\RFA$ and $\C$, 
we can also study differential equations in $\RFA$ where the solutions are mappings from $\RFA$ to $\RFA$. For example, one could consider the following fuzzy differential equation:  
\begin{equation}\label{sistemalinear}
   \left\lbrace \begin{array}{l}
    w'\oplus b\odot w=f(w), \\
   w(z_0)=w_0   
   \end{array} \right.
\end{equation}
where $b\in\RFA$, $w:\RFA\to\RFA$ and $f:\RFA\to\RFA$. 
From Proposition \ref{regrasderivada} and Theorem \ref{teoIntegral}, we obtain that the solution of \eqref{sistemalinear} is given by
\begin{equation}\label{eq:sol_fuzzymapping}
   w(z)=e^{-b\odot (z\ominus z_0)}\left(w_0\oplus  \int^{z}_{z_0}{e^{b\odot (z\ominus z_0)}\odot f(z)\, dz}\right).
\end{equation}
Similar to the complex case, one can expect that the solutions of such differential equations in $\RFA$ are often given in terms of power series.  However, we will leave this study for future work. 
In this work, we only focus on the case in which the solution of a linear equation is given by a curve in $\RFA$.

\section{Application to Two-dimensional System of Fuzzy Differential Equations}

This section presents an initial study on two-dimensional systems of fuzzy differential equations whose main goal is to show some important aspects of the calculus developed in this article. The first one is the chain rule for functions with fuzzy inputs, which involves the product operation between fuzzy numbers. The other one is the division between fuzzy numbers presented in this article.
To this end, we will investigate the dynamics of two-dimensional models through their phase portrait.

Consider the following autonomous system of fuzzy differential equations in terms of the deterministic independent variable $t$: %hale, elstein2005mathematical,wilson): 
\begin{equation}\label{eq:general_bidimension}
\left\{
    \begin{array}{rcl}
    \frac{dx}{dt} &=& g(x(t),y(t))\\
    \frac{dy}{dt} &=& h(x(t),y(t)) 
    \end{array},\right.
\end{equation}
where $g$ and $h$ are fuzzy functions from $\RFA^2$ to  $\RFA$.
This is the case of interaction between two species (see e.g. \cite{bassanezi1988,edelstein2005mathematical,hale2009ordinary}) and its fuzzy case can be found in \cite{de2017first,laiate2021bidimensional} and the references therein.

A qualitative study on the system given in Eq. \eqref{eq:general_bidimension} is the phase portrait. 
In fact, in systems of two-dimensional differential equations with state variables $x$ and $y$, the phase portrait (or plane) illustrates the set of all solutions $x(t)$ and $y(t)$ of the system of equations in question, omitting the independent variable $t$. 
Thus, it is very common to represent the curves of the solutions in the phase space under the form $y=f(x)$, by dividing one equation by another and  considering the chain rule. For example, in the system above, we have 
\begin{equation*}
\frac{dy}{dx} = \frac{h(x,y)}{g(x,y)}.
\end{equation*}
Note that the meaning of $\frac{dy}{dx}$ is given in Definition \ref{def:derivative}, since, in this case, $y$ is regarded as a function from $\RFA$ to $\RFA$. 
Moreover, the quotient $\frac{h(x,y)}{g(x,y)}$ 
exists whenever $g(x,y) \neq 0$ (see Eq. \eqref{eq:division} of Remark \ref{rm:division}). 
%At least two things stand out. One is the meaning of $dy/dx$ without the chain rule and with fuzzy arguments. The other is the guarantee of the quotient $h(x(t),y(t))/g(x(t),y(t))$, which is guaranteed by formula (8) of Remark 1;

Initially, let us illustrate this methodology in the following simple linear system: 
\begin{equation}\label{eq:simples}
\left\lbrace \begin{array}{rl}
  \frac{dx}{dt}  & = - y  \\
  \frac{dy}{dt}  & =  x  \\
\end{array}
\right.
\end{equation}
where $x(t) = r(t) + p(t)A \in \RFA$ and $y(t) = s(t) + q(t)A \in \RFA$. 
From the rule chain and Eq. \eqref{eq:division}, we have 
\begin{equation*}
\frac{dy}{dx} = -\frac{x}{y} \Rightarrow y\odot y'=-x
\Rightarrow \frac{y^2}{2}=-\frac{x^2}{2} \oplus k
\Rightarrow {y^2}\oplus {x^2}=2k,
\end{equation*}
for some  $k \in \RFA$.
So, we have 
\begin{equation*}
\begin{array}{rl}
{y^2}\oplus {x^2}=2k & \Leftrightarrow 
 {(s+qA)^2}\oplus {(r+pA)^2}=2k_{1}+2k_{2}A \\
%     & \Rightarrow{r^2+s^2}+r+pA)^2=2k_{1}+k_{2}A \\
&  \Leftrightarrow ({r^2+s^2})-({p^2+q^2})+2(rp+sq)A= 2k_{1}+2k_{2}A %\\
%& \Rightarrow({r^2+s^2})-({p^2+q^2})= 2k_{1} + 2(k_{2}-sq -rp)A    
\end{array}
\end{equation*}
%y^2}+{x^2}=2k  
%\Rightarrow {(s+qA)^2}+{(r+pA)^2}=2(k_{1}+k_{2}A)
%\Rightarrow{r^2+s^2}+r+pA)^2}=2(k_{1}+k_{2}A
%\Rightarrow({r^2+s^2})-({p^2+q^2})+2(rp+sq)A= 2(k_{1}+k_{2}A
%\Rightarrow({r^2+s^2})-({p^2+q^2})= 2((k_{1}-rp)+k_{2}-sq))A
which produces the following system:
\begin{equation}\label{eq:simples_b}
\left\lbrace \begin{array}{rl}
  ({r^2+s^2})-({p^2+q^2})  & = 2k_{1}  \\
  (rp+sq) = & k_{2}   
\end{array}
\right.
\end{equation}
%\begin{equation*}

Note that if $p=q=0$, that is, the case where there is no uncertainty, the solutions of \eqref{eq:simples_b} are circumferences with center at $(0,0)$ and radius $\sqrt{2k_1}$, $k_1\geq 0$. 
%Note that $\frac{}{}$ (equilibrium point) $(0,0)$. 
If $p \neq 0$ or $q\neq 0$, the family of solutions are circumferences ``disturbed'' by the term $({p^2+q^2})$.

\begin{figure}[H]
\center
\includegraphics[width=0.55\textwidth]{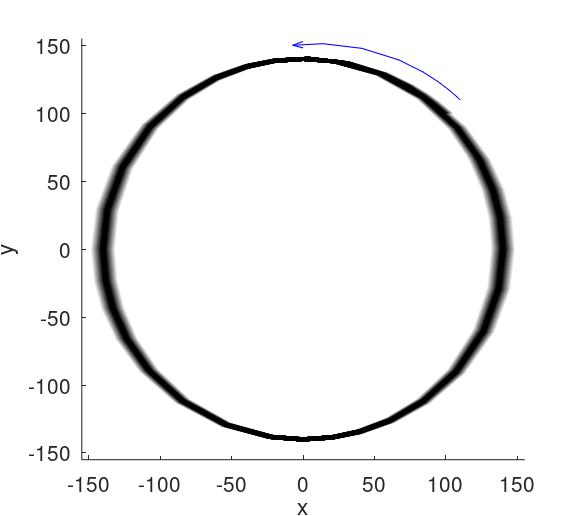}
\caption{Graphical representation of the phase portrait for the pairs given by $\{(x(t),s(t))\}$, where $x(t)$ and $y(t) = s(t) + q(t)A$, $t\in [0,50]$, correspond to the solutions of the system \eqref{eq:simples} with $A = (-1;0;1.01)$ and $x(0) = y(0) = 100 + 2A$. Gray lines, from white to black, represent the endpoints of the $\alpha$-level of $x$ for $\alpha$ varying from 0 to 1, respectively.}
\label{fig:x_simples_phase}
\end{figure}

\begin{figure}[H]
\center
\includegraphics[width=0.55\textwidth]{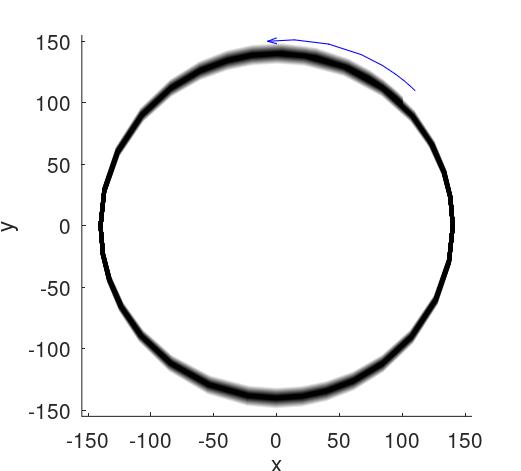}
\caption{Graphical representation of the phase portrait for the pairs given by $\{(r(t),y(t))\}$, where $x(t) = r(t) + p(t)A$ and $y(t)$, $t\in [0,50]$,  correspond to the solutions of the system \eqref{eq:simples} with $A = (-1;0;1.01)$ and $x(0) = y(0) = 100 + 2A$. Gray lines, from white to black, represent the endpoints of the $\alpha$-level of $y$ for $\alpha$ varying from 0 to 1, respectively.}
\label{fig:y_simples_phase}
\end{figure}

Figures \ref{fig:x_simples_phase}  and \ref{fig:y_simples_phase} 
represent the phase portrait of system \eqref{eq:simples} for the cases of $p = 0$ and $q = 0$ with $A = (-1;0;1.01)$ and $x(0) = y(0) = 100 + 2A$.

\begin{figure}[H]
\center
\includegraphics[width=0.55\textwidth]{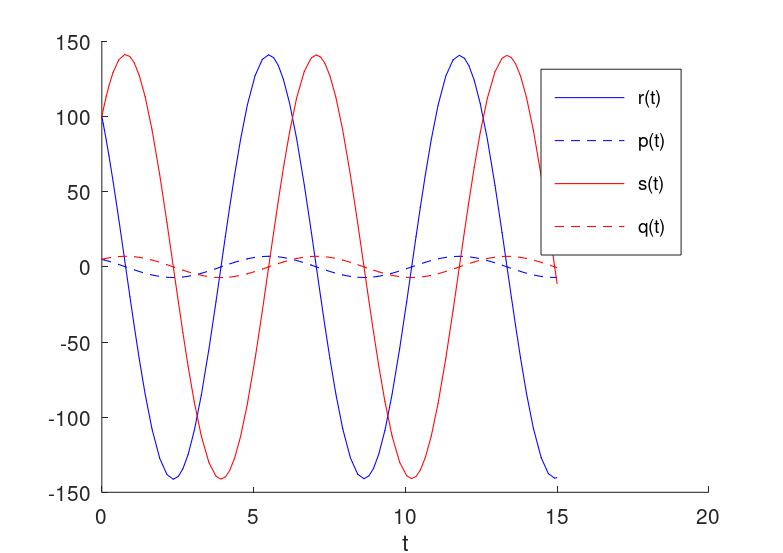}
\caption{Graphical representation of the functions $r,p,s,q$ that are associate with the solution of \eqref{eq:simples} with $A = (-1;0;1.01)$ and $x(0) = y(0) = 100 + 2A$, for $t \in [0,50]$.}
\label{fig:rspq_simples}
\end{figure}

Figure \ref{fig:rspq_simples} presents the the evolution over time of functions $r,p,s,q$ associated with the solutions $x(t) = r(t) + p(t)A$ and $y(t) = s(t) + q(t)A$, for $t \in [0,50]$,  of  Eq. \eqref{eq:simples} with $A = (-1;0;1.01)$ and $x(0) = y(0) = 100 + 2A$. In this case, $r$ and $s$ can be interpreted as 
the trend of the solutions $x$ and $y$, whereas 
$r$ and $s$ as the fuzzy variance/fluctuation coefficient around the corresponding trends. 

\begin{figure}[H]
\center
\includegraphics[width=0.55\textwidth]{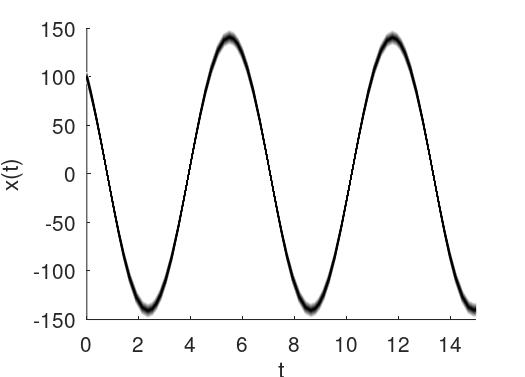}
\caption{Graphical representation of the solution $x$ of \eqref{eq:simples} with $A = (-1;0;1.01)$ and $x(0) = y(0) = 100 + 2A$, for $t \in [0,50]$. Gray lines, from white to black, represent the endpoints of the $\alpha$-level of $x$ for $\alpha$ varying from 0 to 1, respectively.}
\label{fig:x_simples}
\end{figure}

\begin{figure}[H]
\center
\includegraphics[width=0.55\textwidth]{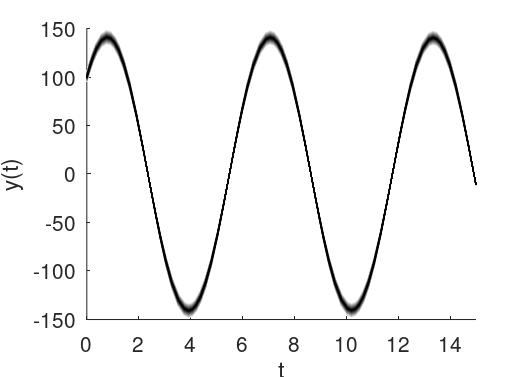}
\caption{Graphical representation of the solution $y$ of \eqref{eq:simples} with $A = (-1;0;1.01)$ and $x(0) = y(0) = 100 + 2A$, for $t \in [0,50]$. Gray lines, from white to black, represent the endpoints of the $\alpha$-level of $y$ for $\alpha$ varying from 0 to 1, respectively.}
\label{fig:y_simples}
\end{figure}

Figures \ref{fig:x_simples}--\ref{fig:y_simples} exhibit the evolution over time of the solutions $x(t) = r(t) + p(t)A$ and $y(t) = s(t) + q(t)A$, for $t \in [0,50]$,  of  the system \eqref{eq:simples} with $A = (-1;0;1.01)$ and $x(0) = y(0) = 100 + 2A$. 

As we can observe in Figures \ref{fig:rspq_simples}--\ref{fig:y_simples}, the solutions $x$ and $y$ of Eq. \eqref{eq:simples} have an oscillatory behaviour. 
In fact, this follows from the fact that positive values of $x$ produce a positive variation in the value of $y$. Similarly, 
due to the negative sign on the left side of the first equation, 
$x$ decreases if $y$ is positive. 
Oscillatory dynamics of two non negative state variables are typical of general prey-predator models. 
In the following subsection, we study a classical and important model in Ecology called Lokta-Volterra which plays a central role in prey-predator systems \cite{bassanezi1988,edelstein2005mathematical}.
%Due to the negative sign in the first equation in the bidimensional system \eqref{eq:simples} has characteristic that $x$ favorece $y$ e este desfavorece $x$. Also, it presents an interesting oscillatory behaviour both in solutions $x(t), y(t)$ and phase portrait. These observations é tipical of predator-prey model general. 
%
%In the follow subsection we will study a very important model in Ecology.

\subsection{The Fuzzy Lotka-Volterra Model}

The Lotka-Volterra equation is a well-known model in Ecology that is used to describe the dynamic of a predator-prey system. 
Here, we consider the generalization of this model by allowing uncertain parameters and variables to take values in the space $\RFA$:
\begin{equation}\label{eq:lotka_volterra}
\left\lbrace \begin{array}{rl}
  x'  & = \alpha \odot x \ominus a \odot x \odot y,  \\
  y'  & = -\beta \odot y \oplus b \odot x \odot y,  \\
  x(0) & = x_0, \\
  y(0) & = y_0, 
\end{array}
\right.
\end{equation}
where $x(t) = r(t) + p(t)A \in \RFA$ and $y(t) = s(t) + q(t)A \in \RFA$ correspond respectively to the population of the prey and predator for $t \in [0,T]$, $\alpha = r_\alpha + q_\alpha A \in \RFA$ is the fuzzy growth rate of the prey, $\beta  = r_\beta + q_\beta A\in \RFA$ is the fuzzy death rate of the predator, $a = r_a + q_a A\in \RFA$ is the predation rate per predator in a unit of time, $\beta = r_b + q_b A \in \RFA$ is the fuzzy growth rate of predator caused by the presence of prey in a unit of time, $x_0 = r_0 + p_0A$ and $y_0 = s_0 + q_0A$ are the fuzzy initial conditions.

Using the fact that every element of $\RFA \setminus \{0\}$ has multiplicative inverse with respect the  multiplication $\odot$, we have that equilibrium points of the system are 
$P_1 = (0,0) \in \RFA^2$ and $P_2 = \left(\frac{\beta}{b},\frac{\alpha}{a}\right) \in \RFA^2$, since 
\[
x'= x \odot (\alpha  \ominus a \odot y) = 0 \Leftrightarrow x = 0 \mbox{ or } y = \frac{\alpha}{a},
\]
and 
\[
y'= y \odot (-\beta  \oplus b \odot x) = 0 \Leftrightarrow y = 0 \mbox{ or } x = \frac{\beta}{b}.
\]
The case of interest is the non trivial equilibrium point $P_2$ which involves divisions of elements of $\RFA$. 

It is worth noting that, in the classical case, the Lokta-Volterra model can be simplified by considering its linearized version around the non trivial equilibrium point. 
This simplification provides an easier system of equations such that its solution in the phase portrait is given by a convenient formula of  an ellipse. 
So, if we consider fuzzy version of its linearization around the non trivial equilibrium point, then we obtain the following system:
%So, before delving deeper into the Lotka-Volterra system \eqref{eq:lotka_volterra}, let us first to consider the following fuzzy version of its linearization around to the non trivial equilibrium point:
\begin{equation}\label{eq:lotka_volterra_linear}
\left\lbrace \begin{array}{rl}
  u'  & = -\frac{a \odot \beta}{b} \odot v,  \\
  v'  & = \frac{b \odot \alpha}{a} \odot u.  
\end{array}
\right..
\end{equation}
The analysis of the phase portrait of this equation is similar to that done for Equation \eqref{eq:simples}. 
In fact, Eq. \eqref{eq:lotka_volterra_linear} corresponds to Eq.  \eqref{eq:simples} if 
$1 = \frac{a \odot \beta}{b}$ and $1 = \frac{b \odot \alpha}{a}$. 
So, we do not delve deeper into this fuzzy linearized system. 
%$c = \frac{a \odot \beta}{b}$ and $d = \frac{b \odot \alpha}{a}$. 
%As can observe in Figures \ref{fig:x_simples_phase} and \ref{fig:y_simples_phase}, as expected from the classical case, its phase portrait resembles a fuzzy ellipse. 
%Therefore, the fuzzy linearized system has been treated above. 

\begin{figure}[H]
\center
\includegraphics[width=0.55\textwidth]{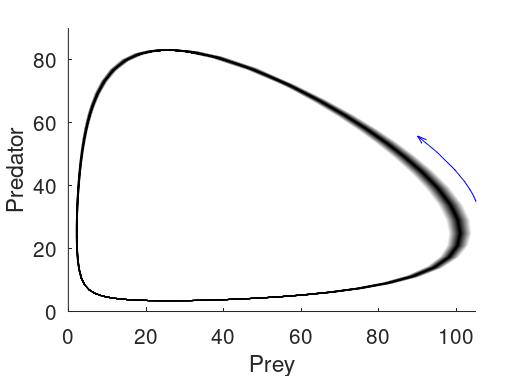}
\caption{Graphical representation of the phase portrait for the pairs given by $\{(x(t),s(t))\}$, where $x(t)$ and $y(t) = s(t) + q(t)A$, $t\in [0,50]$,  correspond to the solutions of the system \eqref{eq:lotka_volterra} with $A = (-0.5;0;0.51)$, $\alpha = 0.25 + 0.001A$, $\beta = 0.18 + 0.003A$, $a = 0.01$, $b=0.007$, $x(0) = x_0 = 100 + 5A$, and $y(0) = y_0 = 30 + 2A$. Gray lines, from white to black, represent the endpoints of the $\alpha$-level of $x$ for $\alpha$ varying from 0 to 1, respectively.}
\label{fig:x_phase_lotka_volterra}
\end{figure}

\begin{figure}[H]
\center
\includegraphics[width=0.55\textwidth]{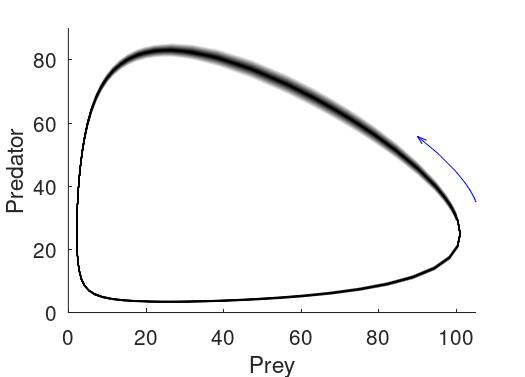}
\caption{Graphical representation of the phase portrait for the pairs given by $\{(r(t),y(t))\}$, where $x(t) = r(t) + p(t)A$ and $y(t)$, $t\in [0,50]$,  correspond to the solutions of the system \eqref{eq:lotka_volterra} with $A = (-0.5;0;0.51)$, $\alpha = 0.25 + 0.001A$, $\beta = 0.18 + 0.003A$, $a = 0.01$, $b=0.007$, $x(0) = x_0 = 100 + 5A$, and $y(0) = y_0 = 30 + 2A$. Gray lines, from white to black, represent the endpoints of the $\alpha$-level of $y$ for $\alpha$ varying from 0 to 1, respectively.}
\label{fig:y_phase_lotka_volterra}
\end{figure}

Suppose that initial conditions of Eq. \eqref{eq:lotka_volterra} are not in equilibrium points. In this case, as we shall see in the phase portrait, the solutions have periodic behavior around to the non-trivial equilibrium point $P_2 = \left(\frac{\beta }{b},\frac{\alpha}{a}\right)$ of Eq. \eqref{eq:lotka_volterra}. 
To analyze the phase portrait, we consider the following equation that follows by applying the rule chain and Eq. \eqref{eq:division}:
\begin{eqnarray*}
y' =\frac{dy}{dx} & = & -\frac{\beta \odot y \oplus b \odot x \odot y}{\alpha \odot x \ominus a \odot x \odot y} = \frac{y \odot \left(-\beta \oplus b \odot x \right)}{x \odot \left( \alpha \ominus a \odot y \right)}  \\
\Leftrightarrow \left(\frac{\alpha \ominus a \odot y}{y}\right)\odot y' & = & \frac{-\beta \oplus b \odot x}{x},
\end{eqnarray*}
%\begin{equation*}
%\frac{dy}{dx} = -\frac{\beta \odot x \oplus b \odot x \odot y}{\alpha \odot x \ominus a \odot x \odot y} \Rightarrow (\alpha \odot x \ominus a \odot x \odot y\odot)y'=-(\beta \odot x \oplus b \odot x \odot y)
%\end{equation*}
by integrating  with respect to $x$ on both sides, the solution of equality above is given implicitly by 
\begin{equation}\label{eq:implicita}
\alpha\odot \mathrm{ln} y \ominus a\odot y =-\beta\odot \mathrm{ln} x \oplus b\odot x \oplus k,
\end{equation}
for some  $k \in \RFA$. Note that the derivation of the left side involves the use of 
chain rule given in Proposition \ref{prop:chain_rfa}.
%To check the last formula above, simply differentiate each side of the equation, according to the derivation rules in $\RFA$.
The parameter $k$ can be determined by the initial conditions. More precisely, we have $k = \alpha\odot \mathrm{ln} y_0 \ominus a\odot y_0 \oplus \beta\odot \mathrm{ln} x_0 \ominus b\odot x_0$ where $x(0) = x_0$ and $y(0) = y_0$ are the the initial conditions.  
Although this system more sophisticated than that of Eq. \eqref{eq:simples}, where the phase portrait is given by Eq. \eqref{eq:simples_b}, the phase portrait of the non-linear Eq. \eqref{eq:lotka_volterra} can be obtained from Eq. \eqref{eq:implicita} by written all parameters and states variables in the form $m+nA \in \RFA$ with $m,n \in \R.$  
Figures \ref{fig:x_phase_lotka_volterra} and \ref{fig:y_phase_lotka_volterra} illustrate some fuzzy trajectories in the phase portrait for the solution of Eq. \eqref{eq:lotka_volterra} with $A = (-0.5;0;0.51)$, $\alpha = 0.25 + 0.001A$, $\beta = 0.18 + 0.003A$, $a = 0.01$, $b=0.007$, and the initial conditions $x(0) = x_0 = 100 + 5A$ and 
$y(0) = y_0 = 30 + 2A$. 

%From Eq. \eqref{eq:implicita}, we can obtain the phase portrait by considering the following auxiliary functions:
%\[
%z(x) = -\beta\odot \mathrm{ln}x \oplus b\odot x \oplus k
%\quad \mbox{ and } \qquad z(y) = \alpha\odot \mathrm{ln} y \ominus a\odot y.
%\]
%
%Let us consider Eq. \eqref{eq:lotka_volterra} with $A = (-0.5;0;0.51)$, $\alpha = 0.25 + 0.001A$, $\beta = 0.18 + 0.003A$, $a = %0.01$, $b=0.007$, and the initial conditions $x(0) = x_0 = 100 + 5A$ and 
%$y(0) = y_0 = 30 + 2A$. 
%With these parameters, we obtain the non trivial equilibrium point 
%\[ 
%P_2 = \left(25.714 + 0.4286A, \; 25 + 0.1A \right) = \left( \;(25.5;25.714;25.933), \; 
%(24.95;25;25.051) \; \right).
%\]
%Figures \ref{fig:x_phase_lotka_volterra} and \ref{fig:y_phase_lotka_volterra} illustrate the phase portrait 
%for the solution of Eq. \eqref{eq:lotka_volterra}.

\begin{figure}[H]
\center
\includegraphics[width=0.6\textwidth]{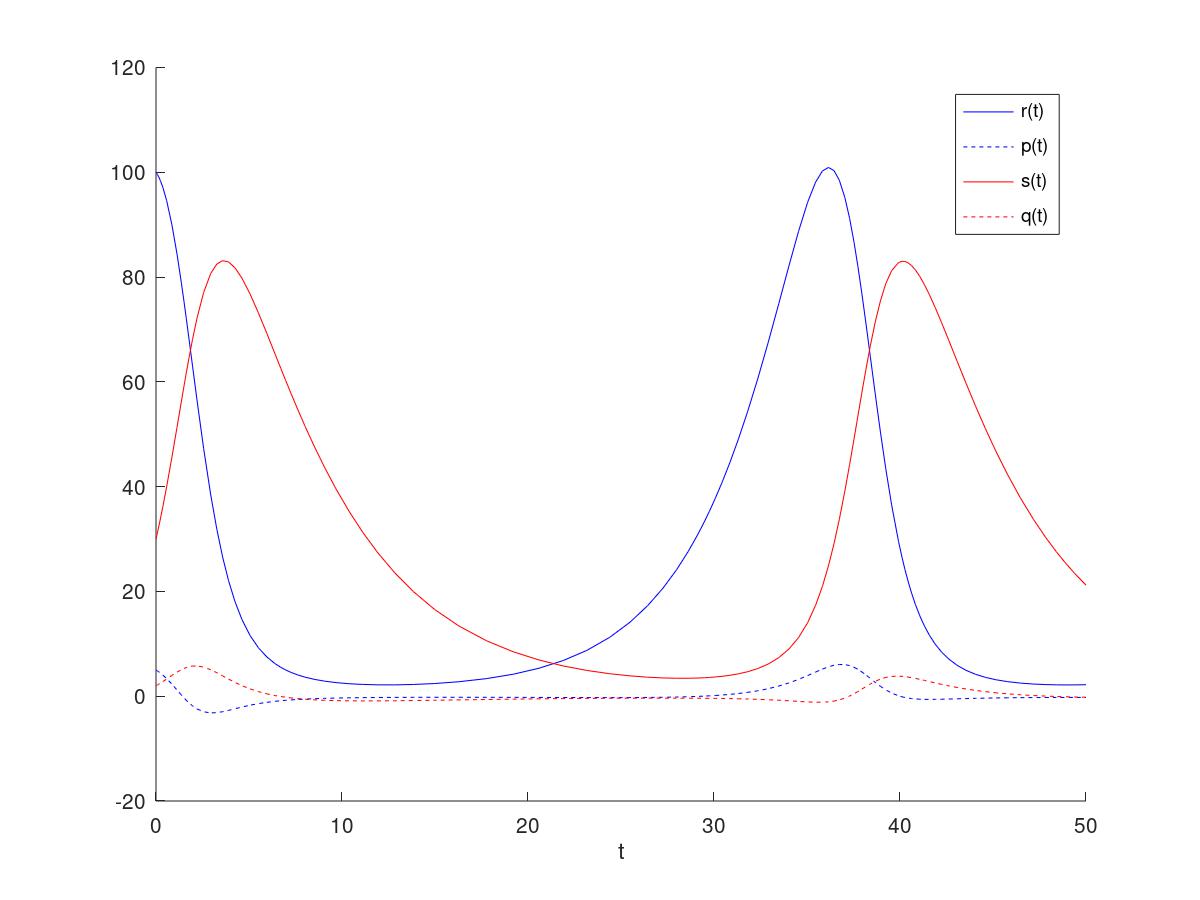}
\caption{Graphical representation of the functions $r,p,s,q$ that are associate with the solution of \eqref{eq:lotka_volterra}  with $A = (-0.5;0;0.51)$, $\alpha = 0.25 + 0.001A$, $\beta = 0.18 + 0.003A$, $a = 0.01$, $b=0.007$, and initial conditions $x(0) = x_0 = 100 + 5A$ and $y(0) = y_0 = 30 + 2A$, for $t \in [0,50]$.}
\label{fig:pqrs_lotka_volterra}
\end{figure}

\begin{figure}[H]
\center
\includegraphics[width=0.55\textwidth]{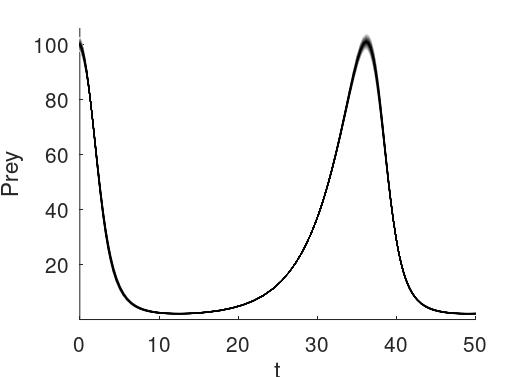}
\caption{Graphical representation of the solution $x$ of \eqref{eq:lotka_volterra} with $A = (-0.5;0;0.51)$, $\alpha = 0.25 + 0.001A$, $\beta = 0.18 + 0.003A$, $a = 0.01$, $b=0.007$, and initial conditions $x(0) = x_0 = 100 + 5A$ and $y(0) = y_0 = 30 + 2A$, for $t \in [0,50]$. Gray lines, from white to black, represent the endpoints of the $\alpha$-level of $x$ for $\alpha$ varying from 0 to 1, respectively.}
\label{fig:x_lotka_volterra}
\end{figure}

\begin{figure}[H]
\center
\includegraphics[width=0.55\textwidth]{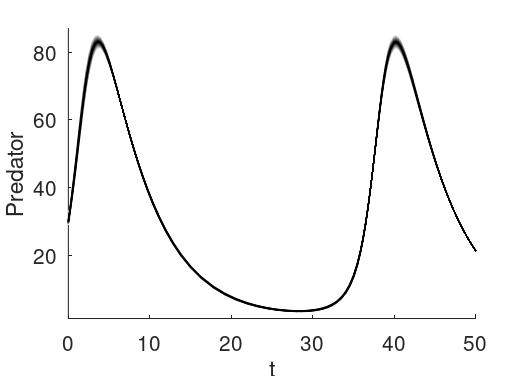}
\caption{Graphical representation of the solution $y$ of \eqref{eq:lotka_volterra} with $A = (-0.5;0;0.51)$, $\alpha = 0.25 + 0.001A$, $\beta = 0.18 + 0.003A$, $a = 0.01$, $b=0.007$, and initial conditions $x(0) = x_0 = 100 + 5A$ and $y(0) = y_0 = 30 + 2A$, for $t \in [0,50]$. Gray lines, from white to black, represent the endpoints of the $\alpha$-level of $y$ for $\alpha$ varying from 0 to 1, respectively.}
\label{fig:y_lotka_volterra}
\end{figure}

By considering $A = (-0.5;0;0.51)$, $\alpha = 0.25 + 0.001A$, $\beta = 0.18 + 0.003A$, $a = 0.01$, $b=0.007$ in 
Eq. \eqref{eq:lotka_volterra} with initial conditions $x(0) = x_0 = 100 + 5A$ and $y(0) = y_0 = 30 + 2A$, we obtain the non trivial equilibrium point 
\[ 
P_2 = \left(25.714 + 0.4286A, \; 25 + 0.1A \right) = \left( \;(25.5;25.714;25.933), \; 
(24.95;25;25.051) \; \right).
\]
Moreover, the evolution of the solutions $x$ and $y$ over time can be obtained by applying the bijection between $\RFA$ and $\C$, that is, by considering $\RFA$  as a fuzzification of $\R$. 
In fact, it follows that Eq. \eqref{eq:lotka_volterra} can be associated with the following classical 
system of ordinary differential equation:
\begin{equation}\label{eq:lotka_volterra_b}
\left\lbrace \begin{array}{rl}
  r'  & = r(r_\alpha - r_as + q_aq) + p(-q_\alpha + r_aq + q_as),  \\
  s'  & = s(-r_\beta + r_br - q_bp) + q(q_\beta - r_bp - q_br), \\
  p'  & = p(r_\alpha - r_as + q_aq) + r(q_\alpha - r_aq - q_as), \\
  q'  & = q(-r_\beta + r_br - q_bp)  + s(-q_\beta + r_bp + q_br), \\
  r(0) & = r_0, \\
  s(0) & = s_0, \\
  p(0) & = p_0, \\
  q(0) & = q_0.
\end{array}
\right..
\end{equation} 
The graphs of the functions $r,s,p,q$ are depicted in Figure \ref{fig:pqrs_lotka_volterra}. 
Moreover, the solutions $x = r + pA$ and $y = s + qA$ are represented in Figures \ref{fig:x_lotka_volterra} and \ref{fig:y_lotka_volterra}, respectively.

In \cite{laiate2021bidimensional}, Laiate {\it et al.}  studied the Lotka-Volterra model from the point of view of the realification of $\RFA$, that is, regarding $\RFA$ as 2-dimension vector space over $\R$. 
In contrast, our approach deals with the Lotka-Volterra model from the point of view of the 
fuzzification of $\R$, that is, regarding $\RFA$ as a field over itself. In this case,  the integral and differential calculus is based on the calculus for analytical functions. 
Moreover, different from what is done in \cite{laiate2021bidimensional}, here we present a brief analysis of the phase portrait of the solution which is only possible due to the chain rule for functions with fuzzy arguments.

\section{Conclusions}\label{Conclusion}

This article presents a fuzzy calculus theory for functions that map fuzzy numbers to fuzzy numbers.  One of the main benefits of this approach is that it allows for modeling phenomena where both inputs and outputs are uncertain. It is worth noting that, in fuzzy differential equations, the inputs of a direction field represent the state of the solution and, therefore, an uncertain quantity.
%To the best of our knowledge, this is the first article that proposes a fuzzy differential and integral calculus for functions where the inputs are fuzzy. 

Although a broader class of fuzzy sets can be considered, this study focuses on the space $\RFA$ that is composed of all fuzzy numbers that are linearly correlated with an asymmetric fuzzy number $A$ (i.e., $\RFA = \{r+qA\mid (r,q) \in \R^2\}$. This space has two important mathematical characteristics: 1) it is a Banach space, which is an ideal framework for developing a theory of differential and integral calculus and, at the same time, 2) it allows us to model phenomena such that uncertainty (or noise) cannot be disregarded in the task of obtaining an adequate solution for the problem at hand.  

Since each $B \in \RFA$ can be written as $B = r+qA$, in analogy to the language used in statistics (or even in stochastic calculus),  $B$ could be understood as the ``composition'' of the real part ``$r$'' and the uncertain part ``$qA$''. In this context, the elements of $\RFA$ could be interpreted from an epistemic point of view. 
On the other hand, from a more technical point of view, such a representation of $B$ also resembles the representation of a complex number, where $A$ plays the role of an imaginary unit ``$i$''.

Thus, the proposed approach is to associate the space $\RFA$ with the set of complex numbers through the bijection $\Phi$, thereby equipping $\RFA$ with the structure of a vector space of complex numbers, $\C$. Consequently, other concepts such as differentiability and integrability are also inherited from $\C$. Additionally, from a practical standpoint, this approach enables the application of conformal maps in the fuzzy domain.

One of the advantages of the proposed approach is to equip $\RFA$ with the multiplication $\odot$ which is induced by the multiplication of complex numbers. Unlike what happens with the $\Psi$-cross product, every element of $\RFA$ that is different from zero has a reciprocal, that is, a multiplicative inverse. Such a property is essential to define the derivative of a function $f:\RFA \to \RFA$ at some $T\in \RFA$ in terms of the limit of $H^{-1}\odot(f(T+H) \ominus f(T))$ as $H\to 0$, where $H^{-1}$ stands for the reciprocal of $H$ with respect to the $H$, i.e, $H\odot H^{-1} =1$.

The vector space $\RFA$ (and the vector space $\C$) over the real field ($\R$) has 2 dimensions, and $\{1,A\}$ (and $\{1,i\}$) is a basis.  However, the set of complex numbers is a vector space over the complex field ($\C$) that has 1 dimension, with its canonical basis being $\{1\}$.
%The vector space $\RFA$ (and the vector space $\C$) over the real field ($\R$) has $2$-dimension and $\{1,A\}$ (and $\{1,i\}$) is a basis. But, the set of complex numbers is a vector space over the complex field ($\C$) that has $1$-dimension and its canonical basis is $\{1\}$. 
Using the isomorphism $\Phi$ between $\C$ and $\RFA$, we show in Theorem \ref{thm:rfa_field} that  %we also obtain that 
$\RFA$ is a field and, hence, it is a $1$-dimension vector space over the filed $\RFA$ with $\{1\}$ as a basis.  
Note that in this last case, the multiplication operation on $\RFA$ is defined as in Equation \eqref{multiplicacao}. 
In contrast, the notion of multiplication is not unique if we regard $\RFA$ as vector space over $\R$  and, therefore, we need to indicate which one is being used.  
%In this case, the solutions of differential equations are obtained  in $\RFA$ 
Depending on the multiplication adopted, the $\Psi$-cross product or that induced by complex numbers, the solutions are qualitatively different. 
For instance, in the linear model with the $\Psi$-cross product, the uncertainty (or the noise) in the growth rate does not affect the dynamic of the real part $x(t)$ (see Equation \eqref{psisol}), whilst, with the multiplication $\odot$ given in  Equation \eqref{multiplicacao}, the uncertainty in this parameter affects both the real and fuzzy parts of the solution, as can be seen in Equation  \eqref{soljuntas}. 
In addition, if we regard $\RFA$ as a $1$-dimension vector space over itself, we can consider differential equations whose solutions are fuzzy mappings, which is the case of the 
solution \eqref{eq:sol_fuzzymapping}. 

We would like to present some brief comments comparing our proposal with the theory of fuzzy calculus based on the concept of generalized Hukuhara derivative or, for short, $gH$-derivative.
On the one hand, at first glance, the notion of the generalized Hukuhara derivative seems to be more general, as it is not restricted to special classes of fuzzy numbers \cite{bede13}.
However, to use it  (mainly when solving fuzzy differential equations), firstly, we need to identify {\it a priori} the switching points of the function, that is, the points where there is a change in the behavior of the widths of the $\alpha$-levels. Secondly, since its calculation is done in the endpoints of $\alpha$-levels, it is necessary to check whether the conditions of Stacking Theorem \cite{bede13} are satisfied to ensure that the derivative is well-defined at a given point. In general, checking such conditions is a very hard task in practice. Therefore, although there is no restriction in the codomain of the functions, its use could be restrictive due to the comments above.         
On the other hand, our approach to fuzzy calculus can be applied to functions taking values in subclasses of fuzzy numbers. These subclasses are subsets of fuzzy numbers whose elements are autocorrelated. Among them, we could highlight those based on linear (or complete) correlation \cite{carlsson2004additions}, for which the space $\RFA$ is included \cite{de2021differential}. 
Under certain conditions, these spaces have the structure of Banach spaces and, therefore, they are closed with respect to the addition, subtraction, and scalar product. This fact eliminates the need to examine the requirements mentioned above when using the $gH$-derivative. Indeed, the theory of differential and integral calculus in Banach space is already well-established. Finally, our approach has explored the Lotka-Volterra model from the perspective of fuzzification of $\R$.
 This involves treating 
$\RFA$ as a field over itself and has included a brief analysis of the solution's phase portrait, facilitated by the chain rule applicable to functions with fuzzy arguments. As a future perspective, this introductory study can be extended for higher dimension of fuzzy differential equation systems and analysis of stability.

\section*{Acknowledgements}
This work was partially supported by CNPq under grants no. 314885/2021-8 and  311976/2023-9, and by FAPESP under grants no. 2022/00196-1 and 2020/09838-0.

\bibliography{references}   %references, veio no artigo:mybibfile

\end{document}